%% file: Unified-Calculus.tex
   \documentclass[twoside,reqno,11pt]{fcaa-var} %
   \input fcaa_style_DG %
   \usepackage{hyperref}
    \def\theequation{\arabic{section}.\arabic{equation}}

 \def\themycopyrightfootnote{\vspace*{3pt}
 \copyright \, 2021\,
   \par  \noindent 
   \hfill  \vspace*{-36pt}
  }

 \title[TOWARDS A UNIFIED THEORY OF FRACTIONAL \dots]
       {TOWARDS A UNIFIED THEORY OF FRACTIONAL \\ [3pt] AND NONLOCAL VECTOR CALCULUS}
 \author[\normalsize M. D'Elia, M. Gulian, H. Olson, G.E. Karniadakis]
        {\normalsize Marta D'Elia $^1$, Mamikon Gulian $^2$, \vspace*{3pt} \break
                     Hayley Olson$^3$, George Em Karniadakis $^4$}

\begin{document}

 \vbox to 1.0cm { \vfill }

 \bigskip \medskip

 \begin{abstract}

Nonlocal and fractional-order models capture effects that classical partial differential equations cannot describe; for this reason, they are suitable for a broad class of engineering and scientific applications that feature multiscale or anomalous behavior.
This has driven a desire for a vector calculus that includes nonlocal and fractional gradient, divergence and Laplacian type operators, as well as tools such as Green's identities,  to model  subsurface transport, turbulence, and conservation laws.
In the literature, several independent definitions and theories of nonlocal and fractional vector calculus have been put forward. Some have been studied rigorously and in depth, while others have been introduced ad-hoc for specific applications.
The goal of this work is to provide foundations for a unified vector calculus by (1) consolidating fractional vector calculus as a special case of nonlocal vector calculus, (2) relating unweighted and weighted Laplacian operators by introducing an equivalence kernel, and (3) proving a form of Green's identity to unify the corresponding variational frameworks for the resulting nonlocal volume-constrained problems.
The proposed framework goes beyond the analysis of nonlocal equations by supporting new model discovery, establishing theory and interpretation for a broad class of operators, and providing useful analogues of standard tools from the classical vector calculus.

 \medskip

{\it MSC 2010\/}: Primary 34B10; Secondary 35R11, 35B40, 26B12

 \smallskip

{\it Key Words and Phrases}:  nonlocal vector calculus; fractional vector calculus; asymptotic behavior of operators

 \end{abstract}

 \maketitle

 \vspace*{-24pt}


\section{Introduction}\label{sec:introduction}

\setcounter{section}{1}
\setcounter{equation}{0}\setcounter{theorem}{0}

The use of nonlocal models in place of their classical differential counterparts has been steadily increasing thanks to their ability to capture effects that classical partial differential equations (PDEs) cannot describe. These effects include multiscale behavior and anomalous behavior such as super- and sub-diffusion, making nonlocal models suitable for a broad class of engineering and scientific applications. Such applications include subsurface transport \cite{Benson2000,katiyar2019general,katiyar2014peridynamic,Schumer2003,Schumer2001}, image processing \cite{Buades2010,Gilboa2007,Lou2010}, multiscale and multiphysics systems \cite{Alali2012,Askari2008,DET19}, magnetohydrodynamic \cite{Schekochihin2008},
phase transitions \cite{Bates1999,Chen_nonlocalmodels,dayal2006kinetics},
finance \cite{Scalas2000,Sabatelli2002}, stochastic processes \cite{Burch2014,DElia2017,Meerschaert2012,Metzler2000,Metzler2004}, and, more recently, fractional backpropagation in training neural networks \cite{Wei2020}.

Nonlocal models are characterized by integral operators acting on the values of a function on {\it nonlocal neighborhoods}; these are usually Euclidean balls of radius $\delta > 0$, which is referred to as the {\it horizon} or {\it interaction radius}. The latter determines the extent of the nonlocal interactions and spans from values much smaller than the size of the domain, including zero in the classical limit, to values much larger than the size of the domain, including infinity in certain fractional-order models. Properties of the kernel defining the integral operator may also reduce the regularity requirements of solutions compared to solutions of classical models, allowing for discontinuities and rough behavior such as crack formation while capturing long-range forces.

Well-known challenges arise from modeling and simulation of nonlocal problems, including the prescription of nonlocal analogues of boundary conditions, i.e., volume or exterior conditions \cite{Cortazar2008,DEliaNeumann2019,DEliaBC2021,Lischke2020}; the unresolved treatment of nonlocal interfaces \cite{Alali2015,Capodaglio2019}; the uncertainty and sparsity of model parameters and data \cite{Antil2019FracControl,burkovska2020,DElia2014DistControl,DElia2016ParamControl,Gulian2019,Pang2019fPINNs,Pang2017discovery,You2020Regression}, and the increasing computational cost as the extent of the nonlocal interactions increases \cite{AinsworthGlusa2017,AinsworthGlusa2018,Aulisa2021,Capodaglio2020DD,DElia-ACTA-2020,DEliaFEM2020,Pasetto2019,Silling2005meshfree,Wang2010,Xu2021FETI}.
In this work, we focus on synthesizing various theoretical frameworks for nonlocal and fractional calculus. Several formulations of nonlocal vector calculus models have been proposed independently, including various fractional-order models. Similarities are evident between certain models, e.g., between the nonlocal gradient operator of Ref.~\cite{Shieh2015} and the fractional gradient introduced in Ref.~\cite{Meerschaert2006}, but comprehensive and rigorous comparisons are lacking, making it difficult to choose the most appropriate fractional model for a specific application.

A formative development of nonlocal calculus for scalar functions was carried out in Ref.~\cite{Gunzburger2010}. Ref.~\cite{Du2013} extended that calculus to vector functions, rigorously derived nonlocal conservation/balance laws, and recast strong forms of nonlocal problems into variational forms, allowing the analysis of nonlocal diffusion equations in a similar framework as their local PDE counterparts. In Ref.~\cite{Du2012}, the authors analyzed nonlocal diffusion equations and made connections with certain types of nonlocal operators (e.g. the graph Laplacian).
More recently, a class of nonlocal operators, originally referred to as {\it weighted} operators \cite{Du2013}, have received increasing attention due to the one-point structure of the weighted nonlocal gradient and its application in mechanics \cite{Du2018Dirichlet} and fluid dynamics \cite{Du2020SPH} contexts. Of particular interest is the stability of the nonlocal problems associated with these operators. Refs. \cite{Du2018Dirichlet,Du2020SPH} showed that radially symmetric and suitably singular kernels are sufficient to guarantee stability, while Ref. \cite{Lee2020} extended this analysis to non-spherical interaction neighborhoods, allowing to relax the singularity requirement. Relevant nonlocal calculus developments include the analysis of nonsymmetric kernels in the context of convection \cite{DElia2017,Tian2017convection}, nonlocal Dirichlet forms \cite{Felsinger2015}, and nonlocal-in-time dynamical systems \cite{Du2017time}, and the extension to nonlocal operators of classical properties of vector field such as Helmholtz decompositions \cite{DElia2020Helmholtz,Lee2020}.
{While this general nonlocal calculus is comparatively recent, fractional-order differential operators are almost as old as their integer-order counterparts \cite{Gorenflo1997,Samko1993}. As a result, the development of fractional calculus literature has often been independent of nonlocal calculus literature; proposals for fractional-order vector calculus include Refs.~\cite{Meerschaert2006,silhavy2020fractional,Tarasov2008}. Not only is there a lack of rigorous comparison between these theories, but there is also a lack of discussion of fractional vector calculus as it may relate to the aforementioned nonlocal vector calculus of Refs.~\cite{Du2013} and \cite{Gunzburger2010}.
Thus, for example, fractional-order Green's identities have been reported in Ref.~\cite{Dipierro2017}, but no connection has been drawn to Green's identities in the nonlocal vector calculus literature, such as as in Ref.~\cite{Du2013}.

Our aim is to unify fractional-order and nonlocal vector calculus, and develop tools within this unified framework to study well-posedness and other properties.
We now outline our paper. In Section \ref{sec:preliminaries}, we review necessary aspects of nonlocal and fractional vector calculus. In Section \ref{sec:fractional-nonlocal}, we prove the equivalence of Cartesian and directional fractional vector calculus and show they are special cases of general weighted nonlocal vector operators for infinite-range interactions. In the same vein, we relate unweighted nonlocal and fractional Green's identity, prove equivalence of the composition of fractional divergence and gradient with the fractional Laplacian, and present results on the asymptotic behavior of truncated fractional operators, which is relevant for numerical simulations.
In Section \ref{sec:unweighted-weighted-equivalence}, we continue the effort initiated in Refs. \cite{Du2018Dirichlet,Du2020SPH,Lee2020,mengesha2014peridynamic} towards the development of a unified characterization of nonlocal Laplacian operators and formally prove the equivalence of weighted and unweighted nonlocal Laplacian operators by introducing an {\it equivalence kernel}, for which we discuss general properties and show equivalence to well-known fractional kernels. In Section \ref{sec:unified-variational}, we first prove a weighted nonlocal Green identity. Then, we exploit this result to unify the variational framework for weighted operators with the framework for unweighted operators and provide well-posedness results for a class of singular kernels, for which fractional operators are a special instance. This further confirms previous stability results for truncated nonlocal operators \cite{Du2018Dirichlet,Du2020SPH}. Finally, in Section \ref{sec:Tarasov}, we review other instances of fractional vector calculus, focusing in particular on the theories put forth by Tarasov \cite{Tarasov2008} and {\v{S}}ilhav{\`y} \cite{silhavy2020fractional}, comparing them to the unified framework proposed here.

The resulting theory provides a new framework for the discovery of nonlocal models that admit clear interpretation, properties, and are well-posed by construction; see Refs.~\cite{Pang2020,You2020Regression,You2021} for recent examples. The established connections between nonlocal and fractional operators guide the choice of specific fractional models and their discretizations in the context of numerical simulations, see, e.g.,  \cite{DiLeoni-2020} where a new fractional gradient for turbulence modeling is proposed and for which our theory provides insights on its asymptotic behavior. This foundational work invites future studies into classes of kernels and function spaces for which a unified theory could be established with full rigor for specific applications.}

\vspace*{-1pt} 
\section{Preliminaries on nonlocal operators}\label{sec:preliminaries}

\setcounter{section}{2}
\setcounter{equation}{0}\setcounter{theorem}{0}

In this section we describe unweighted and weighted nonlocal operators as well as fractional operators, reviewing results that will be relevant for the rest of the article. The bulk of our new results concerns weighted nonlocal vector calculus, but we review the unweighted nonlocal vector calculus in some detail in order to establish a relation between weighted and unweighted problems in Section \ref{sec:unified-variational}.

We note that several integrals may have to be intended in the principal value (p.v.) sense; however, for the sake of notation, we avoid adding the p.v. symbol in front of the integral sign.

\vspace*{-8pt}
\subsection{Interaction Radius}

Key to all the definitions in the section will be the antisymmetric vector-valued function $\bm{\alpha}(\mathbf{x}, \mathbf{y})$, known as the interaction kernel, which is used to describe the nonlocal interaction between points $\mathbf{x}$ and $\mathbf{y}$.
Throughout, we write the interaction kernel as
\begin{equation}\label{eq:alpha-rho}
\bm{\alpha}(\xb,\yb) = \bm{\rho}(\xb,\yb)\mathds{1}(|\xb-\yb|\leq\delta).
\end{equation}
We explicitly use $\delta$ to quantify the ``range of interaction'', and the antisymmetric vector-valued function $\bm{\rho}$ to describe the interaction within this range. Thus, there is zero interaction when $\xb$ and $\yb$ are further than $\delta$ apart.

\begin{figure}[t]
\centering
\includegraphics[width=0.7\textwidth]{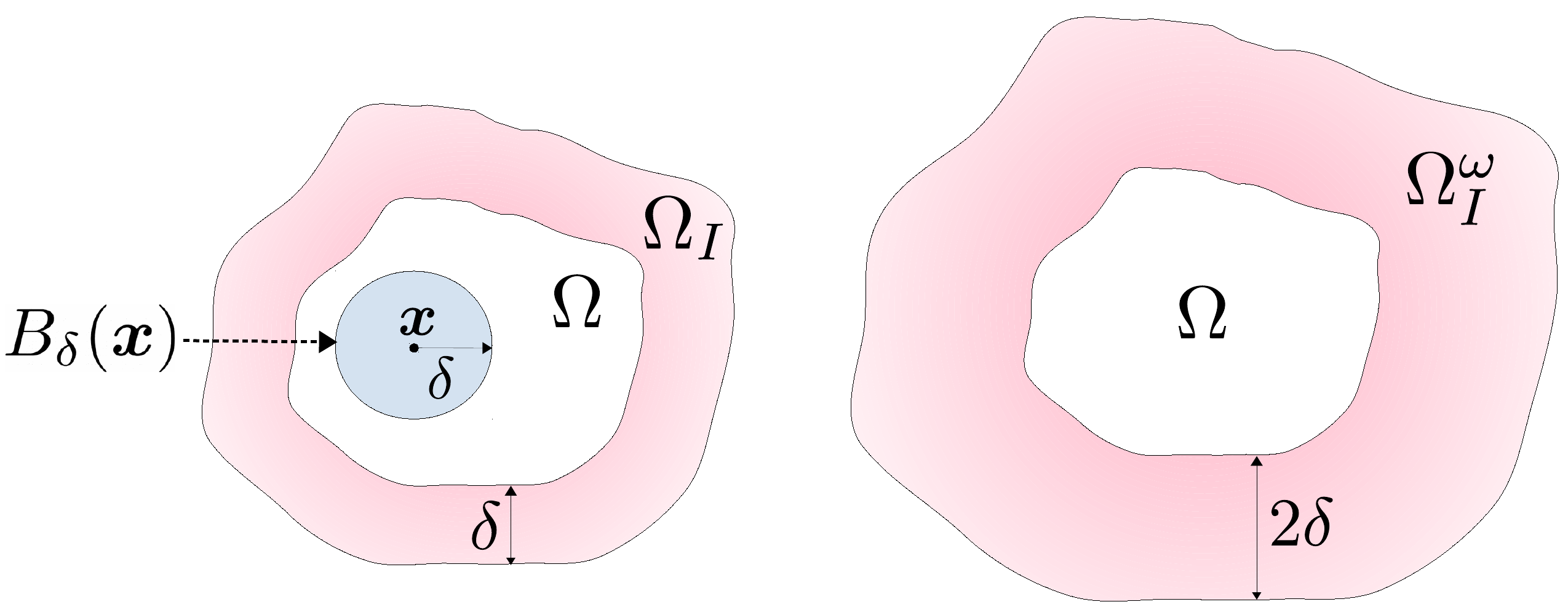} 
\\[-0.3em]
\caption{Illustration and comparison of the unweighted interaction domain $\Omega_I$ (\emph{left}) and weighted interaction domain $\Omega_I^\omega$ (\emph{right}). The ball $B_\delta(\xb)$ represents the support of the kernel $\alphab$. The interaction radius $\delta$ determines the extent of the nonlocal interactions. The domain $\Omega$ is the same in both illustrations. The configuration on the left is used for \emph{unweighted} nonlocal Poisson problems, while the thicker interaction domain on right is required for \emph{weighted} nonlocal Poisson problems.}
\label{fig:domains}
\end{figure}

We refer to operators defined with $\delta < \infty$ as \emph{truncated} operators, and operators with $\delta = \infty$ as \emph{untruncated} operators.
Tightly connected to the function $\bm{\alpha}$ is the {\it interaction domain} $\Omega_I$, which consists of points outside a domain of interest $\Omega$ that interact with points inside $\Omega$:
\begin{equation*}
\omgi = \{ \yb \in \mbRn \setminus\omg \ \text{such that} \
\xb \; \hbox{interacts with} \; \yb\ \text{for some}\ \xb \in \omg\}.
\end{equation*}
For now, this definition is purposefully vague; we shall see below that the precise definition of $\Omega_I$ differs between the unweighted and weighted cases.
In Figure \ref{fig:domains} we report two possible configurations of a domain $\Omega$ and its interaction domain.

\vspace*{-6pt}
\subsection{Unweighted operators}\label{sec:unweighted_operators} 

Although unweighted nonlocal operators were introduced in early works by Gilboa and Osher (see, e.g., Ref.~\cite{Gilboa2008}), we refer to Refs.~\cite{Du2012,Du2013,Gunzburger2010} as the standard references providing a rigorous nonlocal vector calculus, and we follow their notation and nomenclature.

Let $\alphab:\mathbb{R}^n\times\mathbb{R}^n\to\mathbb{R}^n$, for $n=1,2,3$, be an anti-symmetric two-point vector function defined as in \eqref{eq:alpha-rho}. For $\vb:\mathbb{R}^n\times\mathbb{R}^n\to\mathbb{R}^n$, the nonlocal {\it unweighted divergence} $\mathcal{D}\vb:\mathbb{R}^n\to\mathbb{R}$ is defined as
\vskip -10pt %
\begin{equation}\label{eq:unw-div}
\begin{aligned}
\mathcal{D}\vb(\xb) := \int_{\mathbb{R}^n} (\vb(\xb,\yb)+\vb(\yb,\xb))\cdot\alphab(\xb,\yb)d\yb. \hspace{1cm}
\end{aligned}
\end{equation}
For $u:\mathbb{R}^n\to\mathbb{R}$ the nonlocal {\it unweighted gradient}, $\mathcal G u:\mathbb{R}^n\times\mathbb{R}^n\to\mathbb{R}^n$, the negative adjoint$^{i}$ \footnote{$^{i}$ This corresponds to -$\mathcal{D}^*u(\xb,\yb)$ in the notation of  Ref.~\cite{Du2013}; see equation (3.14a) of that article. We use the notation $\mathcal{G}$ rather than $-\mathcal{D}^*$ to establish notational similarity between classical vector calculus and nonlocal vector calculus. Note that Ref.~\cite{Du2013} introduces an operator $\mathcal{G}$ in the context of unweighted nonlocal vector calculus which is entirely different from equation \eqref{eq:unw-grad} above.} of \eqref{eq:unw-div}, is defined as
\begin{equation}\label{eq:unw-grad}
\mathcal{G}u(\xb,\yb) = (u(\yb)-u(\xb))\alphab(\xb,\yb). \hspace{1cm}
\end{equation}
For the kernel $\gamma=\alphab\!\cdot\!\alphab$, we define the nonlocal {\it unweighted Laplacian} as the composition of unweighted nonlocal divergence and gradient, i.e.
\begin{equation}\label{eq:unw-lapl}
\mcL u(\xb) = \mcD\mcG u(\xb) =
2\int_\mbRn (u(\yb)-u(\xb))\gamma(\xb,\yb)d\yb.
\end{equation}
Note that, by construction, $\gamma$ in the equation above is always positive.

\subsubsection{A nonlocal unweighted vector calculus.} 
The operators \eqref{eq:unw-div}, \eqref{eq:unw-grad} and \eqref{eq:unw-lapl} have been rigorously studied in Ref.~\cite{Du2013}, where a nonlocal vector calculus mimicking the local counterpart was developed.
This theory includes a nonlocal unweighted Green identity and the strong and weak forms of an unweighted truncated nonlocal Poisson problem. Since we wish to develop a similar framework for weighted operators, we review these results.

For posing a problem of the form $-\mcL u = f$ in a bounded domain $\Omega$, it is necessary to define an \emph{interaction domain}. Considering the definition \eqref{eq:unw-lapl} and the support of the interaction kernel \eqref{eq:alpha-rho}, to evaluate $\mcL u(\xb)$ for $\xb \in \Omega$ it is necessary to use values $u(\xb)$ for $\xb \in \Omega_\delta$, where
\begin{equation}
\Omega_\delta =
\{\yb \in \mathbb{R}^n \text{ such that } |\xb - \yb| \le \delta
\text{ for some } \xb \in \Omega
\}.
\end{equation}
The interaction domain associated with the unweighted nonlocal Laplacian in \eqref{eq:unw-lapl} is defined as
\begin{align}\label{eq:interaction-dom-unweighted}
\Omega_I &=
\{\yb\in\mbRn\!\setminus\!\Omega
\text{ such that }
\alpha(\xb,\yb)\neq 0 \text{ for some } \xb\in\Omega\} \\
&=
\omg_\delta \setminus \omg.
\end{align}

We introduce the strong form of a such an unweighted nonlocal Poisson problem.
For $f:\Omega\to\mbR$ and $g:\Omega_I\to\mbR$, find $u$ such that
\begin{equation}\label{eq:bound-truncated-prob}
\left\{\begin{aligned}
-\mcL u (\xb) = f (\xb) & \quad \xb\in\Omega\\
u(\xb) = g(\xb)        & \quad \xb\in\Omega_I,
\end{aligned}\right.
\end{equation}
where the second condition in \eqref{eq:bound-truncated-prob} is the nonlocal counterpart of a Dirichlet boundary condition for PDEs; as such, we refer to it as {\it Dirichlet volume constraint} $^{ii}$. \footnote{$^{ii}$ Neumann and Robin volume constraint could also be considered, but, for the sake of simplicity, we only treat the homogeneous Dirichlet case. We refer to Ref. \cite{Du2012} for definition and analysis of the Neumann problem and to, e.g., Ref. \cite{DEliaNeumann2019} for its numerical treatment.}
Here, we see that the set $\omgi$ plays the role of nonlocal boundary.
In fact, Ref.~\cite{Du2012} shows that for problems of the form \eqref{eq:bound-truncated-prob}, the volume constraint on the solution in $\Omega_I$ is required to guarantee well-posedness.

For simplicity and without loss of generality we analyze the homogeneous case $g=0$ on $\Omega_I$; note that the results below can be extended to the non-homogeneous case using ``lifting'' arguments (see, e.g., Ref.~\cite{DElia2014DistControl}).
As for PDEs, a fundamental result for recasting the strong form \eqref{eq:bound-truncated-prob} into a variational problem is nonlocal integration by parts. Ref.~\cite{Du2013} proves the following nonlocal first Green identity:
\vskip -10pt%
\begin{multline} \label{eq:unweighted-Green}
  \int\limits_\Omega -\mcL u(\xb) \, v(\xb) \,d\xb
\\
= \int\limits_\omgomgi\int\limits_\omgomgi \mcG u(\xb,\yb) \cdot \mcG v(\xb,\yb) \,d\yb\,d\xb
+ \int\limits_\omgi \mcD(\mcG u)(\xb)\,v(\yb) \,d\xb.
\end{multline}
Multiplying \eqref{eq:bound-truncated-prob} by a test function $v$ such that $v = 0$ on $\Omega_I$ and integrating over the domain $\Omega$ yields
\vskip-10pt %
\begin{align}\label{eq:bound-truncated-weak}
0 &=  \int\limits_\Omega (-\mcL u -f)\, v \,d\xb
\\
&= \int\limits_\omgomgi\int\limits_\omgomgi \mcG u \cdot \mcG v \,d\yb\,d\xb
+ \int\limits_\omgi \mcD(\mcG u)\,v \,d\xb
- \int\limits_\Omega f\, v \,d\xb.
\end{align}
Thus, the weak form of the nonlocal Poisson problem reads as follows. For $f\in V'(\omgomgi)$ and $g\in V(\omgomgi)|_\omgi$, find $u\in V(\omgomgi)$ such that $u=0$ in $\omgi$ and
\begin{equation}\label{eq:bound-truncated-weak-forms}
\mcA(u,v) = \mcF(v), \;\;\forall\, v\in V(\omgomgi),
\end{equation}
\vskip -3pt \noindent %
where
\vskip -12pt
\begin{equation}\label{eq:A-F-V}
\begin{aligned}
\mcA(u,v) & = \int_\omgomgi\int_\omgomgi \mcG u \cdot \mcG v \,d\yb\,d\xb,\\[2mm]
\mcF(v) & =\int_\Omega f\, v \,d\xb,\\[2mm]
V(\omgomgi) & =\{v\in L^2(\omgomgi): \vertiii{v}<\infty\;
                {\rm and} \; v|_\omgi = 0\}.
\end{aligned}
\end{equation}
Here, the {\it energy norm} $\vertiii{\cdot}$ is defined as
\begin{equation}\label{eq:unweighted-energy}
\vertiii{v}^2 = \int_\omgomgi\int_\omgomgi
|\mcG v|^2\,d\yb\,d\xb.
\end{equation}
The spaces $V'$ and $V|_\omgi$ are, respectively, the dual space and trace space on $\omgi$ of $V$. Note that the bilinear form $\mcA(\cdot,\cdot)$ defines an inner product on $V(\omgomgi)$ and that $\vertiii{u}^2= \mcA(u,u)$. This fact, the continuity of $\mcA$ and $\mcF$ and the Lax-Milgram theorem, yield the well-posedness of the weak form \eqref{eq:bound-truncated-weak-forms}, \cite{Du2012}.

\vspace*{-8pt}
\subsection{Weighted operators} 
\label{sec:intro_weighted_operators}
The operators \eqref{eq:unw-div} and \eqref{eq:unw-grad} are the building blocks of weighted nonlocal operators introduced in Ref.~\cite{Du2013}. The major difference between the unweighted and weighted families of operators is that the weighted gradient yields a one-point function, and the weighted divergence takes one-point functions as its argument. Throughout, we let $\omega:\mathbb{R}^n\times\mathbb{R}^n\to\mathbb{R}$ be a nonnegative, symmetric scalar function known as the \emph{weight}. For $\vb:\mathbb{R}^n\to\mathbb{R}^n$, the nonlocal {\it $\omega$-weighted divergence} $\mcD_\omega\vb:\mathbb{R}^n\to\mathbb{R}$ is defined as
\vskip -10pt %
\begin{align}
\label{eq:w-div}
\mathcal{D}_{\omega} \vb(\xb) :&= \mcD(\omega(\xb,\yb)\vb(\xb))
\\
&=
\int_{\mathbb{R}^n} (\omega(\xb,\yb)\vb(\xb) + \omega(\yb,\xb)\vb(\yb))\cdot\alphab(\xb,\yb)d\yb.
\end{align}
Ref.~\cite{Du2013} shows that an equivalent definition is given by
\vskip -10pt %
\begin{equation}\label{eq:w-div-alternative-def}
\mathcal{D}_{\omega} \vb(\xb) =
\int_{\mathbb{R}^n} (\vb(\yb) - \vb(\xb))\cdot
\alphab(\xb,\yb)\omega(\xb,\yb)d\yb ;
\end{equation}
this will be useful later on in the paper.
For $u:\mathbb{R}^n\to\mathbb{R}$, the nonlocal {\it $\omega$-weighted gradient} $\mcG_\omega u:\mbRn\to\mbRn$ is defined as
\vskip -10pt %
\begin{align}\label{eq:w-grad}
\mcG_\omega u(\xb) :&= \int_\mbRn\mcG u(\xb,\yb) \omega(\xb,\yb)d\yb \\
&=
\int_\mbRn
(u(\yb)-u(\xb))\alphab(\xb,\yb)
\omega(\xb,\yb)d\yb .
\end{align}
Ref.~\cite{Du2013} shows that the latter is the negative adjoint of the former. As done in the unweighted case, we define the nonlocal {\it $\omega$-weighted Laplacian} as the composition of \eqref{eq:w-div} and \eqref{eq:w-grad}, i.e., for $\xb\in\Omega$
\begin{equation}\label{eq:w-lapl}
\begin{aligned}
\mcLw u(\xb) &= \mcD_\omega\mcG_\omega u(\xb)\\
&= \int_\mbRn \Big[\omega(\xb,\yb)\mcG_\omega u(\xb) + \omega(\yb,\xb) \mcG_\omega u(\yb)\Big] \cdot\alphab(\xb,\yb)d\yb\\
&
\begin{multlined}
=\int_\mbRn \left[
\int_\mbRn
(u(\zb)-u(\xb))\alphab(\xb,\zb)
\omega(\xb,\zb)d\zb  \right.\\ +
\left.\int_\mbRn
(u(\zb)-u(\yb))\alphab(\yb,\zb)
\omega(\yb,\zb)d\zb
\right] \cdot\alphab(\xb,\yb)\omega(\xb,\yb)d\yb.
\end{multlined}
\end{aligned}
\end{equation}
Note that although we specified $\xb \in \Omega$ in defining $\mathcal{L}_\omega u(\xb)$, the above expansion of $\mcD_\omega \mcG_\omega u(\xb)$ is valid for all $\xb \in \mathbb{R}^n$. However, as a matter of notation, we refer to $\mcD_\omega \mcG_\omega u(\xb)$ as $\mcL_\omega u(\xb)$ only for $\xb \in \Omega$. This is to draw analogy with classical theorems of vector calculus, for example, when discussing the nonlocal analogue of Green's identity in Section \ref{sec:weighted-Green}. Just as for unweighted operators, we refer to weighted operators with $\delta < \infty$ in \eqref{eq:alpha-rho} as truncated, and those with $\delta = \infty$ as untruncated.

Although we restrict ourselves to the nonlocal operators discussed above, related operators have been studied in the literature.
In Ref.~\cite{Alali2015Calculus}, the authors generalize the nonlocal divergence \eqref{eq:w-div} and show that the calculus of Ref.~\cite{Du2013} is a specific instance of their generalized theory.
In Ref.~\cite{Radu2019} the authors introduce a {\it doubly nonlocal} Laplace operator similar to \eqref{eq:w-lapl} and draw connections between it and the standard Laplacian; in contrast to \eqref{eq:w-lapl}, it is not expressed as a composition of nonlocal divergence and gradient operators.

\smallskip

Given an open bounded set $\Omega\subset\mbRn$
we will consider equations of the form
$-\mcLw u = f$ in $\Omega$. As for the unweighted case, problems involving such equations require a volume constraint for an appropriate interaction domain, which can be read off from \eqref{eq:alpha-rho}
 and \eqref{eq:w-lapl}.
Due to the presence of the double integral in $\mcLw$, the interaction domain associated with the operator in \eqref{eq:w-lapl} is defined as
\begin{align}
\begin{split}
\label{eq:interaction-dom-weighted}
\Omega^{\omega}_{I} &= \omg_{2\delta} \setminus \omg
= \big\{
\yb \in \mbRn \setminus \Omega
\text{ such that }
|\xb-\yb| \leq 2\delta
\text{ for some }
\xb \in \Omega
\big\}.
\end{split}
\end{align}
This interaction domain $\Omega^{\omega}_{I}$  for the weighted case is in contrast to the ``thinner'' interaction domain $\omgi$ defined by \eqref{eq:interaction-dom-unweighted} for the unweighted case; this is illustrated in Figure \ref{fig:domains}.
Although this point is sometimes not clearly stated in the nonlocal calculus literature, it is required to obtain well-posed problems and to develop numerical schemes for nonlocal mechanics as in
Refs. \cite{delia2017handbook} and \cite{seleson2016convergence}, where the interaction domain associated with the {\it peridynamic} truncated weighted operator has thickness $2\delta$.

\subsubsection{Nonlocal weighted Dirichlet problems} 
The weighted operators presented in this subsection have not been analyzed in depth as their unweighted counterparts. In this paragraph we introduce strong and weak forms of a weighted Poisson problem and refer to Section \ref{sec:weak-equivalence} for their analysis. We do not distinguish between truncated and untruncated as the latter case simply corresponds to $\delta=\infty$.
We introduce the strong form of a weighted, nonlocal Poisson problem. For $f:\Omega\to\mbR$ and $g:\Omega^{\omega}_{I}\to\mbR$, find $u$ such that
\begin{equation}\label{eq:weight-bound-truncated-prob}
\left\{\begin{aligned}
-\mcLw u = f & \quad \xb\in\Omega\\
u = g        & \quad \xb\in{\Omega^{\omega}_{I} }
\end{aligned}\right. ,
\end{equation}
where the second condition in \eqref{eq:weight-bound-truncated-prob} is still referred to as Dirichlet volume constraint. Next, by multiplying \eqref{eq:weight-bound-truncated-prob} by a test function  $v$ and integrating over the domain $\Omega$, we have the following weak form of \eqref{eq:weight-bound-truncated-prob}:
\begin{equation}\label{eq:weight-truncated-weak}
  \int_\Omega (-\mcLw u -f)\, v \,d\xb = 0.
\end{equation}
The latter has not been analyzed in full generality due to the lack of a complete study of a nonlocal weighted Green identity, which would allow expression of \eqref{eq:weight-truncated-weak} in terms of an energy principle, as in \eqref{eq:bound-truncated-weak-forms}. We define the {\it weighted energy} as
\begin{equation}\label{eq:weighted-energy}
\vertiii{v}_{\omega}^2 = \int_{\Omega \cup \Omega^{\omega}_{I} } (\mcGw v)^2\,d\xb.
\end{equation}
For now, we do not refer to energy \eqref{eq:weighted-energy} as a norm; conditions on $\alphab$ and $\omega$ will guarantee such property and are investigated in Section \ref{sec:weak-equivalence}. {We point out, however, that some studies in this direction have been conducted in previous works. In Ref.~\cite{Du2018Dirichlet} the authors show the equivalence of \eqref{eq:weighted-energy} and the unweighted energy norm \eqref{eq:unweighted-energy} in the case of periodic problems in one dimension. A weighted nonlocal Green identity relevant to \eqref{eq:weighted-energy} was asserted in Ref.~\cite{du2013analysis} and used to derive Euler-Lagrange equations for a peridynamic theory of linear elasticity, but a proof was not provided. A nonlocal integration by parts result for a type of weighted gradient operator was introduced in Ref. \cite{mengesha2016characterization}. In Ref.~\cite{Shieh2015} the authors pose the equivalence of the weighted \eqref{eq:weighted-energy} and unweighted \eqref{eq:unweighted-energy} energies as an open problem (listed as Open Problem 1.9).} We generalize the result of Ref.~\cite{Du2018Dirichlet} and partially answer Open Problem 1.9 of Ref.~\cite{Shieh2015} in Section \ref{sec:weak-equivalence}.

\vspace*{-6pt}
\subsection{Fractional calculus operators}
\label{sec:intro_fractional_operators}

We now introduce two definitions of fractional order 
operators; in both cases, these operators are one-point functions, as for nonlocal weighted operators. To initially distinguish between the two families of operators, we refer to the first class of operators as \emph{directional}, and the second class as \emph{Cartesian}.

\subsubsection{Directional operators}\label{sec:directional_operators} 
In the vector calculus of Ref.~\cite{Meerschaert2006}, the {\it fractional divergence} of a vector field $\vb(\xb)$ is defined as
\begin{equation}\label{eq:frac-grad}
\mathbbmss{grad}^s u(\xb) = \int_{\thetab = 1} \thetab D_{\thetab}^s u(\xb) M(d\thetab), \quad \xb \in\mathbb{R}^n,
\end{equation}
where $s\in(0,1)$ and $M(\cdot)$ is a measure that can account for anisotropy. The directional derivative $D^{s}_{\thetab} u(\xb)$
is the one-dimensional positive Riemann-Louiville fractional derivative$^{iii}$
\footnote{$^{iii}$
While $D^s_\thetab \vb(\xb)$ is defined in Ref.~\cite{Meerschaert2006} in terms of the Fourier transform, this representation is clearly explained on pages 170 and 175 in Ref. \cite{Meerschaert2012}. Ref. \cite{Meerschaert2006} also discusses representations in terms of fractional integrals.}
in $t$, i.e.
\begin{equation*}
\frac{d^s w}{dt^s}=\frac{1}{\Gamma(1-s)}\frac{d}{dt}
\int_{0}^{\infty} w(t-t') \left(t'\right)^{-s} dt'
\end{equation*}
applied to the function $w: \mathbb{R} \rightarrow \mathbb{R}$, $w(t) = u(\xb + t\thetab)$, for fixed $\xb$ and $\thetab$, and then evaluated at $t=0$:
\vskip -12pt
\begin{equation}
D^{s}_{\thetab} u(\xb) =\frac{d^s w}{dt^s}\Big|_{t = 0}\,.
\end{equation}
Ref.~\cite{Meerschaert2006} points out that
$\mathbbmss{grad}^s u(\xb)$ has Fourier transform
\begin{equation}\label{eq:grad_fourier_transform}
\mcF[\mathbbmss{grad}^s u](\xib)
= \int_{|\thetab| = 1} \thetab(i\xib \cdot \thetab)^s
 \hat{u}(\xib) \; M(d\thetab).
\end{equation}

Similarly to \eqref{eq:frac-grad}, the {\it fractional divergence}
of a vector field $\vb(\xb)$ is defined in Ref.~\cite{Meerschaert2006} as
\vskip -11pt
\begin{equation}\label{eq:frac-div}
{\mathbbmss{div}}^s \vb(\xb) = \int_{\thetab = 1} \thetab \!\cdot\! D_{\thetab}^s \vb(\xb) M(d\thetab), \ \; \xb \in\mathbb{R}^n.
\end{equation}
We understand the directional derivative $D_{\thetab}^s \vb(\xb)$ of the vector field $\vb$ to mean the vector
\vskip -12pt
\begin{equation}
D_{\thetab}^s \vb(\xb) =
\big(D_{\thetab}v_1(\xb), \hdots, D_{\thetab}v_n(\xb)\big).
\end{equation}
Ref.~\cite{Meerschaert2006} points out that
$\mathbbmss{div}^s u(\xb)$ has Fourier transform
\begin{equation}\label{eq:div_fourier_transform}
\mcF[\mathbbmss{div}^s \vb](\xib)
= \int_{|\thetab| = 1} (i\xib \cdot \thetab)^s
\hat{\vb}(\xib) \cdot \thetab \; M(d\thetab).
\end{equation}

\begin{remark}
In this subsection we have kept the measure $M(\cdot)$ to remain consistent with Ref.~\cite{Meerschaert2006} and to emphasize the generality of the directional fractional vector calculus. However, in the remainder of the paper we only consider the case of uniform measure
\begin{equation}\label{eq:uniform_M}
M(d\bm{\theta}) = d\bm{\theta}
\end{equation}
to allow for comparisons to other theories of fractional vector calulus. Thus, when we refer to equations \eqref{eq:frac-grad}, \eqref{eq:grad_fourier_transform}, \eqref{eq:frac-div} and \eqref{eq:div_fourier_transform} throughout the article, we assume that \eqref{eq:uniform_M} holds.
\end{remark}

\subsubsection{Cartesian operators} \label{sec:cartesian_operators} 
Another family of fractional vector calculus operators has been the focus of works such as Refs.~\cite{Mazowiecka2018,Ponce2016,Shieh2015,Shieh2017,silhavy2020fractional}. We refer to these operators as \emph{Cartesian} operators, and introduce them below. It should be noted that none of the listed references establish or assert a relationship between the Cartesian operators and the directional operators introduced in the previous subsection. This will be the focus of Section \ref{sec:fractional-nonlocal-equivalence}. We define
\vskip -12pt %
\begin{align}\label{eq:frac-Cart}
\begin{split}
{\rm grad}^s  u(\xb) &=
\int_{\mbRn}\left[ u(\xb) - u(\yb) \right]\,
\frac{\xb-\yb}{|\xb-\yb|}
\frac{1}{|\xb-\yb|^{n+s}}
d\yb ,
\\
{\rm div}^s \vb(\xb) &=
\int_{\mbRn}
\left[ \vb(\xb) - \vb(\yb) \right]\,
\frac{\xb-\yb}{|\xb-\yb|}
\frac{1}{|\xb-\yb|^{n+s}}
d\yb .
\end{split}
\end{align}
Note that it is possible to introduce {\it truncated} fractional operators by multiplying the integrands in \eqref{eq:frac-Cart} by the indicator function of $|\xb-\yb|$ over the ball of radius $\delta$. We refer to such operators as {\it truncated fractional divergence and gradient} and denote them by $\Dsd$ and $\Gsd$:
\begin{align}
\begin{split}
\label{eq:truncated-frac-Cart}
{\rm grad}^s_\delta  u(\xb) &=
\int_{\mathbb{R}^n}\left[ u(\xb) - u(\yb) \right]
\frac{\xb-\yb}{|\xb-\yb|}
\frac{\mathds{1}(|\xb-\yb|\le\delta)}{|\xb-\yb|^{n+s}}
d\yb ,
\\
\Dsd \vb(\xb) &=
\int_{\mathbb{R}^n} \left[ \vb(\xb) - \vb(\yb) \right]
\cdot\frac{\xb-\yb}{|\xb-\yb|}
\frac{\mathds{1}(|\xb-\yb|\le\delta)}{|\xb-\yb|^{n+s}}
d\yb.
\end{split}
\end{align}

The (Riesz) fractional Laplacian, which appears in several of our results below, is an example of a nonlocal Laplacian operator, as discussed in Section \ref{sec:frac-special-case}. It is defined \cite{cai2019,Lischke2020} by
\vskip -11pt
\begin{equation}
\label{eq:riesz_Laplacian_rn}
(-\Delta)^s u = C_{n,s}
\int_{\mathbb{R}^n} \frac{u(\xb) - u(\yb)}{|\xb-\yb|^{n+2s}}d\yb,
\end{equation}
\vskip -3pt \noindent
where
\vskip -13pt
\begin{equation}\label{eq:frac_laplc_const}
C_{n,s}=
\frac{4^{s} \Gamma\left(s+\frac{n}{2}\right)}
{\pi^{n/2}|\Gamma(-s)|}.
\end{equation}

\subsubsection{Tempered operators} 
Tempered fractional operators were introduced in Ref. \cite{Sabzikar2015} and have been studied in the context of stochastic jump processes. In this context, tempering refers to multiplying the stable L\'evy distribution associated to a fractional operator, which feature algebraic power-law tails, by an exponentially decaying factor; analytically, this corresponds to multiplying the kernel of the fractional operator by such a factor. The \emph{tempered fractional Laplacian} is defined as \cite{deng2018boundary,zhang2018riesz}
\vskip -10pt
\begin{equation}
\label{eq:tempered_Laplacian_rn}
(-\Delta)_\lambda^s u = C_{n,s}
\int_{\mathbb{R}^n} \frac{u(\xb) - u(\yb)}{|\xb-\yb|^{n+2s}}
e^{-\lambda|\xb-\yb|}
d\yb,
\end{equation}
where $C_{n,s}$ is the same as \eqref{eq:frac_laplc_const}.

In analogy to \eqref{eq:tempered_Laplacian_rn} and\eqref{eq:riesz_Laplacian_rn},
we introduce definitions of fractional gradient and divergence in Cartesian form: for $\lambda>0$, $\vb:\mathbb{R}^n\to\mathbb{R}^n$, and $u:\mathbb{R}^n\to\mathbb{R}$ the {\it tempered fractional divergence and gradient} are given by
\vskip-11pt %
\begin{align}
\begin{split}
\label{eq:tempered-fractional}
\Gsl u(\xb) &=
\int_{\mbRn}
\left[ u(\xb) - u(\yb) \right]
\frac{\xb-\yb}{|\xb-\yb|}
\frac{e^{-\lambda|\xb-\yb|}}{|\xb-\yb|^{n+s}}
d\yb , \\
\Dsl \vb(\xb) &=
\int_{\mbRn}
\left[ \vb(\xb) - \vb(\yb) \right]
\cdot
\frac{\xb-\yb}{|\xb-\yb|}
\frac{e^{-\lambda|\xb-\yb|}}{|\xb-\yb|^{n+s}}
d\yb .
\end{split}
\end{align}
We study these operators as a specific case of the unified nonlocal vector calculus in Section \ref{sec:gamma-consistency}.

\vspace*{-1pt} 
\subsubsection{Properties of fractional operators.} 
In this paragraph we report several technical results related to the fractional divergence and gradient operators introduced above. These results will be used in the proofs of the main results in the following sections.

We introduce the vector$^{iv}$
\footnote{$^{iv}$ In what follows, vector spaces are denoted by bold symbols.}
fractional Sobolev space ${\bf H}^s(\mbRn)$: we say that a vector field $\vb$ belongs to ${\bf H}^s(\mbRn)$ if each component $\vb_i$ of $\vb$ belongs to $H^s(\mbRn)$, i.e.
$
\int (1+|\xib|^2)^s |\hat{\vb}_i|^2 d\xib< \infty,
$
for each $i = 1, 2,..., n$. Equivalently, for $\hat\vb=\mcF[\vb]$, $\mcF[\cdot]$ being the Fourier transform,
\begin{equation}\label{eq:Hs-componentwise}
\text{$\vb \in {\bf H}^s(\mathbb{R}^n)$ \  if and only if } \
\int (1+|\xib|^2)^s |\hat{\vb} \cdot \hat{\vb}| d\xib < \infty.
\end{equation}

The following lemma provides a result on integrability of the weight $|\xb-\yb|^{-(n+t)}$ away from the singularity $\xb$. It is later used to prove the convergence of the truncated operators \eqref{eq:truncated-frac-Cart} to their untruncated counterparts \eqref{eq:frac-Cart} in Section \ref{sec:convergence}. The proof, based on spherical integration, is omitted.
\vspace*{-3pt}

\begin{lemma}\label{upsilon_lemma}
For $n=1,2,\ldots$ and $t > 0$,
\begin{equation}\label{eq:Upsilon}
\int_{\mbRn \setminus B_\delta(\xb)}
\frac{1}{|\xb-\yb|^{n+t}} d\yb = \frac{\Upsilon_{n,t}}{\delta^t},
\quad \text{where} \quad
\Upsilon_{n,t} =
\frac{\pi^\frac{n}{2}n}{t \Gamma\left(\frac{n}{2}+1\right)}.
\end{equation}
\end{lemma}

The following result is used to prove consistency of
$\mathbbmss{div}^s \mathbbmss{grad}^s$ with
the fractional Laplacian \eqref{eq:riesz_Laplacian_rn} in Theorem \ref{thm:frac-div-grad-lapl}. It follows from an adaptation of the argument of Example 6.24 in Ref.~\cite{Meerschaert2012} and is proven in \ref{sec:integral-equality-constant}. It is also used below to establish mapping properties of the directional forms of fractional gradient and divergence.
\vspace*{-4pt}

\begin{lemma}\label{lem:integral-equality-constant}
There exists a negative constant $D_{n,s}$ such that, for $s\in(0,1)$,
\vskip -11pt%
\begin{equation}\label{eq:integral-equality-constant}
\int_{|\thetab|=1} \int_{|\thetab'|=1} \thetab \cdot \thetab' (i\thetab \cdot \xib)^{s} (i\thetab' \cdot \xib)^{s} d\thetab d\thetab'
 = D_{n,s}|\xib|^{2s}.
\end{equation}
\end{lemma}

The next two lemmas prove important mapping properties (domain and range) of the directional forms of the fractional operators \eqref{eq:frac-grad} and \eqref{eq:frac-div}. The proof of Lemma \ref{gradient_mapping} is based on the same arguments as those of Lemma \ref{divergence_mapping} and is hence omitted.

\vspace*{-3pt}

\begin{lemma}
\label{divergence_mapping}
Let $r \ge 0$.
If $\vb \in {\bf H}^{r+s}(\mbRn)$, then
${\rm div}^{s} \vb \in H^{r}(\mbRn)$.
\end{lemma}

\proof
We show that, for $\vb \in {\bf H}^{r+s}(\mbRn)$, we have
\begin{equation*}
\int |\mathcal{F} \left[ \mathbbmss{div}^{s} \vb \right] (\xib)|^2
(1+|\xib|^2)^{r}
d\xib < \infty.
\end{equation*}
Using the Fourier representation \eqref{eq:div_fourier_transform}, we have
\vskip -11pt%
\begin{align*}
\int |\mathcal{F} \left[ \mathbbmss{div}^{s} \vb \right]& (\xib)|^2 d\xib \\
&\le
 \int
|\hat{\vb}(\xib)|^2 \left| \int \thetab (i\thetab \cdot \xib)^{s} d\thetab \right|^2
(1+|\xib|^2)^{r}
d\xib \\
&=
\int
|\hat{\vb}(\xib)|^2
\left[
\int \int \thetab (i\thetab \cdot \xib)^{s}
\cdot \thetab' (i\thetab' \cdot \xib)^{s} d\thetab d\thetab'
\right]
(1+|\xib|^2)^{r}
d\xib \\
&=
\int
|\hat{\vb}(\xib)|^2 |\xib|^{2s}
(1+|\xib|^2)^{r}
d\xib \\
&\le
\int
|\hat{\vb}(\xib)|^2 (1+|\xib|^2)^{s}
(1+|\xib|^2)^{r}
d\xib \\
&=
\int
|\hat{\vb}(\xib)|^2 (1+|\xib|^2)^{r+s}
d\xib < \infty,
\end{align*}
where the third inequality follows from Lemma \ref{lem:integral-equality-constant} and the last inequality follows from \eqref{eq:Hs-componentwise}.
\proofend

\begin{lemma}
\label{gradient_mapping}
Let $r \ge 0$. If $u \in H^{r+s}(\mbRn)$, then ${\rm grad}^{s} u \in {\bf H}^r(\mbRn)$.
\end{lemma}


\section{Relations between fractional and weighted nonlocal calculus}\label{sec:fractional-nonlocal}

\setcounter{section}{3}
\setcounter{equation}{0}\setcounter{theorem}{0}

In this section we prove several equivalence results. First, we show that the fractional divergence and gradient in directional and Cartesian form are equivalent and that their composition is equivalent to the well-known fractional Laplacian operator. Second, we show that they can be expressed as weighted nonlocal operators. Continuing with the theme of connecting concepts from nonlocal and fractional literature, we also show the equivalence of a recent fractional Green  identity with the standard nonlocal unweighted Green  identity in \eqref{eq:unweighted-Green}. Finally, we show that truncated fractional operators converge to their untruncated counterparts.

\vspace*{-5pt}
\subsection{Equivalence of directional and Cartesian fractional vector calculus} \label{sec:fractional-nonlocal-equivalence}

As a first step towards showing the equivalence of nonlocal and fractional operators, we show that the directional and Cartesian definitions are indeed the same, up to a constant. Then in Section \ref{sec:frac-special-case}, we will show that they are also instances of weighted nonlocal operators.
\vspace*{-3pt}

\begin{theorem}\label{lemma:polar-Cart-equivalence}
For $\vb \in {\bf H}^s(\mathbb{R}^d)$,
$u \in H^s(\mathbb{R}^d)$, and
$M(\thetab)=\thetab$, the directional and Cartesian definitions of fractional divergence and gradient are equivalent, i.e.
\begin{equation}\label{eq:polar-Cart-equivalence}
\begin{aligned}
\mathbbmss{grad}^s  u(\xb) & = G_s {\rm grad}^s  u(\xb),\\
\mathbbmss{div}^s \vb(\xb) & = G_s {\rm div}^s \vb(\xb),
\end{aligned}
\end{equation}
where $G_s={s}/{\Gamma(1-s)}$.
\end{theorem}

\proof
We treat the fractional gradient first. We rewrite the fractional directional derivative in generator form$^{v}$ \footnote{$^{v}$ See Ref.~\cite{Meerschaert2012}, p. 175. Here, the ``generator form of the directional derivative'' refers to the generator form of the positive Riemann-Louiville derivative of $v(t) = u(\xb + t\thetab)$, evaluated at $t = 0$.},
\begin{equation*}
D^{s}_{\thetab} u(\xb)=
\frac{s}{\Gamma(1-s)}\int_0^\infty
[u(\xb) - u(\xb - r\thetab)]r^{-s-1} dr.
\end{equation*}
Then we can write the fractional gradient as
\begin{equation*}
\mathbbmss{grad}^s u(\xb) =
\frac{s}{\Gamma(1-s)}\int_{|\thetab| = 1} \int_0^\infty\thetab
[u(\xb) - u(\xb - r\thetab)]r^{-s-1} dr d\thetab.
\end{equation*}
If we let $\yb = r\thetab$, noting that $d\yb = r^{n-1} dr d\thetab$ , we have
\begin{equation*}
\mathbbmss{grad}^s u(\xb) = \frac{s}{\Gamma(1-s)}
\int_\mbRn (u(\xb)- u(\xb-\yb))
\frac{\yb}{|\yb|}\frac{1}{|\yb|^{n+s}} d\yb.
\end{equation*}
Performing a further change of variables
$\xb - \yb \leftrightarrow \yb$, we have
\begin{align*}
\mathbbmss{grad}^s u(\xb) &=
\frac{s}{\Gamma(1-s)} \int_\mbRn (u(\xb)-u(\yb))
\frac{\xb-\yb}{|\xb-\yb|}\frac{1}{|\xb-\yb|^{n+s}}\,d\yb
\\
&= G_s \,{\rm grad}^s u(\xb).
\end{align*}

Next, we treat the divergence operator. We write the directional derivative of the vector field $\vb$
in the generator form,	
\begin{equation*}
D^{s}_{\thetab} \vb(\xb)
=
\frac{s}{\Gamma(1-s)}
\int_0^\infty
[\vb(\xb) - \vb(\xb - r\thetab)]
r^{-s-1} dr,
\end{equation*}
and insert into the fractional divergence \eqref{eq:frac-div} to get
\begin{equation*}
\mathbbmss{div}^s \vb(\xb) =
\frac{s}{\Gamma(1-s)}
\int_{|\thetab| = 1}
\int_0^\infty
\thetab \cdot
[\vb(\xb) - \vb(\xb - r\thetab)]
r^{-s-1} dr d\thetab.
\end{equation*}
Then, performing the same change of variables
$\yb = r\thetab$ with $d\yb = r^{n-1} dr d\thetab$, we obtain
\vskip -10pt%
\begin{equation*}
\mathbbmss{div}^s  \vb(\xb) =
\frac{s}{\Gamma(1-s)}\int_{\mbRn}
[\vb(\xb) - \vb(\xb - \yb)] \cdot \frac{\yb}{|\yb|}
\frac{1}{|\yb|^{n+s}}d\yb.
\end{equation*}
With the change of variables $\yb \leftrightarrow \xb - \yb$, we have
\begin{align*}
\mathbbmss{div}^s  \vb(\xb) &=
\frac{s}{\Gamma(1-s)}
\int_{\mbRn} (\vb(\xb) - \vb(\yb))\cdot
\frac{\mathbf{x-y}}{|\mathbf{x-y}|}\frac{1}{|\mathbf{x-y}|^{n+s}} d\yb \\
&= G_s\, {\rm div}^s\vb(\xb).
\end{align*}
\endproof

\vspace*{-10pt}

\begin{remark}
The above result implies that $\mathbbmss{div}^s$ and
$\mathbbmss{grad}^s$ are essentially interchangeable with ${\rm div}^s$ and ${\rm grad}^s$, respectively. In the proofs of our results in the remainder of the paper, we freely use whichever representation is most convenient. The constant $G_s$ must be taken into account when changing representations.
\end{remark}

\vspace*{-12pt}
\subsection{Consistency of $\mathbbmss{div}^s \mathbbmss{grad}^s$ with the fractional Laplacian} \label{sec:divs-grads-lapls}

The next theorem proves that the composition of fractional divergence and gradient is equivalent to the fractional Laplacian operator. This extends results in Ref.~\cite{Meerschaert2006}, where compositions with integer-order operators with fractional-order operators were studied.
\vspace*{-3pt}

\begin{theorem}\label{thm:frac-div-grad-lapl}
For $u \in H^{2s}(\mathbb{R}^d)$ the following equivalence holds:
\begin{equation}\label{eq:F-WN-equivalence}
\mathbbmss{div}^s \mathbbmss{grad}^s u = D_{n,s} (-\Delta)^s u,
\end{equation}
where $D_{n,s}$ is the negative real-valued constant in Lemma \ref{lem:integral-equality-constant}.
\end{theorem}

\proof
Lemma \ref{gradient_mapping} implies $\mathbbmss{grad}^s u \in \mathbf{H}^{s}(\mathbb{R}^d)$, and in turn
Lemma \ref{divergence_mapping} implies
\vskip -12pt
\begin{equation}
\mathbbmss{div}^s\mathbbmss{grad}^s u \in L^2(\mathbb{R}^d).
\end{equation}
Given the Fourier transforms of the fractional divergence \eqref{eq:div_fourier_transform} and gradient \eqref{eq:grad_fourier_transform} operators, the composition $\mathbbmss{div}^s \mathbbmss{grad}^s u(\xb)$ has Fourier transform
\vskip -10pt %
\begin{equation*}
\begin{aligned}
    \int_{|\thetab| = 1} (i\xib \cdot \thetab)^s
&   \left[\int_{|\bm{\theta'}| = 1}
    \bm{\theta'}(i\xib \cdot \bm{\theta'})^s
    \hat{u}(\xib) d\bm{\theta'}
    \right] \cdot \thetab d\thetab \\
& = \int_{|\thetab| = 1} \int_{|\bm{\theta'}| = 1}
    (i\xib \cdot \thetab)^s (i\xib \cdot \bm{\theta'})^s
    \hat{u}(\xib)\bm{\theta'} \cdot \thetab d\bm{\theta'}d\thetab \\
& = \hat{u}(\xib)\int_{|\thetab| = 1} \int_{|\bm{\theta'}| = 1}
    (i\xib \cdot \thetab)^s (i\xib \cdot \bm{\theta'})^s
    (\bm{\theta'} \cdot \thetab) d\bm{\theta'}d\thetab\\
& = D_{n,s}|\xib|^{2s}\hat{u}(\xib),
\end{aligned}
\end{equation*}
where $D_{n,s}$ is a negative real-valued constant, by Lemma \ref{lem:integral-equality-constant}.
\endproof

\begin{remark}\label{remark:consistency_L_divgrad}
The same result, up to the constant in \eqref{eq:F-WN-equivalence}, can be also proven using the equivalence between $\mathbbmss{div}^s \mathbbmss{grad}^s$ and ${\rm div}^s {\rm grad}^s$ provided by Theorem \ref{lemma:polar-Cart-equivalence} and the \emph{equivalence kernel} for the latter weighted nonlocal operator; see Section \ref{sec:gamma-consistency}, Theorem \ref{thm:frac-equivalence-kernel}. This consistency between the Cartesian form ${\rm div}^s {\rm grad}^s$ and $-(\Delta)^s$ was also obtained in Refs. \cite{silhavy2020fractional} and \cite{horvath1959some} using Riesz potentials; see Section \ref{sec:silhavy} for a summary. Finally, we point out that a similar result is stated in Ref. \cite{Mazowiecka2018}, where the authors claim that the equivalence holds in a distributional sense, but do not provide a proof nor a reference.
\end{remark}


\vspace*{-14pt}
\subsection{Fractional vector calculus as a special case of weighted nonlocal vector
 calculus}\label{sec:frac-special-case}

The following theorem shows that for a specific choice of kernel function and weight, fractional divergence and gradient are special cases of the nonlocal weighted divergence and gradient operators.
\vspace*{-4pt}

\begin{theorem}\label{thm:fractional_special_case_nonlocal}
Let $\vb \in {\bf H}^s(\mathbb{R}^d)$ and $u \in H^s(\mathbb{R}^d)$.
For the weight function and kernel
\vskip -10pt
\begin{equation}\label{eq:equivalence-w-alpha}
\begin{array}{l}
\displaystyle
\omega= |\yb-\xb| \phi(|\yb-\xb|),
\;\; \mbox{with} \;\;
\phi(|\yb-\xb|) = |\xb-\yb|^{-(n+1+s)}, \\[2mm]
\displaystyle
\rhob(\xb,\yb) = \frac{\yb-\xb}{|\yb-\xb|},
\;\; \mbox{with} \;\;
\delta=\infty,
\end{array}
\end{equation}
the fractional divergence and gradient operators can be identified with the weighted nonlocal operators,
\vskip -10pt
\begin{equation}\label{eq:polar-Cart-weighted-equivalence}
\begin{aligned}
{\rm div}^s \vb(\xb) & = \mcDw \vb(\xb)\\[2mm]
{\rm grad}^s  u(\xb) & = \mcGw u(\xb).
\end{aligned}
\end{equation}
\end{theorem}

\proof
For $\delta=\infty$, we have $B_\delta(\xb)=\mbRn$ for all $\xb\in\mbRn$, and substituting $\omega$ and $\rhob$ in \eqref{eq:equivalence-w-alpha} into $\mcGw$ gives
\vskip -10pt %
\begin{equation*}
\begin{aligned}
\mcGw u(\xb)
& = \int_\mbRn (u(\yb)-u(\xb))\alphab(\xb,\yb)\omega(\xb,\yb)d\yb \\
& = \int_\mbRn (u(\yb)-u(\xb))\frac{\yb-\xb}{|\yb-\xb|}
    \frac{1}{|\yb-\xb|^{n+s}} d\yb \\
& = {\rm grad}^s u(\xb).
\end{aligned}
\end{equation*}
\vskip -3pt \noindent %
Similarly, substitution of $\omega$ and $\rhob$ in \eqref{eq:equivalence-w-alpha} into the definition of the weighted nonlocal divergence in \eqref{eq:w-div-alternative-def} yields
$
\mcDw \vb(\xb)=  {\rm div}^s\vb(\xb).
$
\endproof 

In contrast, the fractional Laplacian arises most naturally from the \emph{unweighted} calculus: in \eqref{eq:unw-div} -- \eqref{eq:unw-lapl}, set $\delta = \infty$, $\Omega_I = \mathbb{R}^n \setminus \Omega$, and
\begin{equation}\label{eq:dipierro_alpha}
\bm{\alpha}(\xb, \yb) = \frac{\mathbf{y} - \mathbf{x}}
{|\mathbf{y} - \mathbf{x}|}
\frac{1}{|\mathbf{x} - \mathbf{y}|^{\frac{n}{2}+s}}.
\end{equation}
Then
\vskip -12pt
\begin{equation}
\gamma(\mathbf{x}, \mathbf{y}) = \frac{1}{|\mathbf{x} - \mathbf{y}|^{n+2s}},
\end{equation}
\vskip -3pt \noindent
and
\vskip -10pt
\begin{equation}
\mathcal{L}u(\mathbf{x}) = -2
\int_{\mathbb{R}^n}
\frac{u(\mathbf{x}) - u(\mathbf{y})}
{|\mathbf{x} - \mathbf{y}|^{n+2s}} \, d\yb
=
\frac{1}{C_{n,s}}(-\Delta)^{s}u(\xb),
\end{equation}
where the integrals are understood in the principal value sense.
As explained in Remark \ref{remark:consistency_L_divgrad}, $(-\Delta)^{s}u$ also arises as the weighted composition $\mcLw = \mcDw \mcGw$ with the interaction kernel given by \eqref{eq:equivalence-w-alpha}.

\vskip 2pt

Although Theorem \ref{thm:fractional_special_case_nonlocal} makes clear that the treatment of fractional-order vector calculus requires the weighted nonlocal framework we point an important example of a result in fractional vector calculus being consistent with nonlocal vector calculus, namely that a recent fractional-order Green identity reported by Ref.~\cite{Dipierro2017} is a special case of the unweighted nonlocal Green identity \eqref{eq:unweighted-Green}.
Specifically, Ref.~\cite{Dipierro2017} proposed a detailed concept of fractional flux and volume constraint, which was implemented studied numerically in Ref.~\cite{Antil2019FracControl}. Ref.~\cite{Dipierro2017} defines, for $\xb \in \mathbb{R}^n \setminus \Omega$,
\begin{equation}\label{eq:dipierro_neumann_MAIN}
\mathcal{N}_s u(\xb) =
\int_{\Omega}
\frac{u(\mathbf{x}) - u(\mathbf{y})}
{|\mathbf{x} - \mathbf{y}|^{n+2s}}
\,d\yb
\end{equation}
as a fractional Neumann derivative. The intuition is that this operator represents interaction between $\xb \in \mathbb{R}^n \setminus \Omega$ and $\yb \in \Omega$, i.e., ``flux'' between $\Omega$ and
$\mathbb{R}^n \setminus \Omega$. The same article reports the fractional Green identity
\vskip -12pt
\begin{multline}\label{eq:dipierro_MAIN}
\int_\Omega (-\Delta)^s u(\xb) \, v(\xb) \,d\xb
+\int_{\mathbb{R}^n\setminus\Omega} \
v(\xb)\,
\mathcal{N}_s u(\xb)
\,d\xb \\
=
\frac{C_{n,s}}{2}
\int \int_{\left(\mathbb{R}^n\right)^2 \setminus \left(\mathbb{R}^n \setminus \Omega\right)^2}
\frac{\left(u(\xb) - u(\yb)\right)\left(v(\xb) - v(\yb)\right)}{|\mathbf{x} - \mathbf{y}|^{n+2s}}
\,d\yb\,d\xb.
\end{multline}
On the other hand, the unweighted Green identity \eqref{eq:unweighted-Green} with \eqref{eq:dipierro_alpha} reads
\vskip-10pt %
\begin{multline}\label{eq:unweighted-Green-fractional_MAIN}
\frac{2}{C_{n,s}}
\int_\Omega (-\Delta)^s u(\xb) \, v(\xb) \,d\xb
=
\int_{\mathbb{R}^n} \int_{\mathbb{R}^n}
\frac{\left(u(\xb) - u(\yb)\right)\left(v(\xb) - v(\yb)\right)}{|\mathbf{x} - \mathbf{y}|^{n+2s}}
\,d\yb\,d\xb \\
-2 \int_{\mathbb{R}^n\setminus\Omega} \
v(\xb)
\int_{\mathbb{R}^n}
\frac{u(\mathbf{x}) - u(\mathbf{y})}
{|\mathbf{x} - \mathbf{y}|^{n+2s}}
\,d\yb\,d\xb.
\end{multline}
In Appendix \ref{sec:appendix_green}, we show that \eqref{eq:unweighted-Green-fractional_MAIN} yields \eqref{eq:dipierro_MAIN}. We mention that nonlocal flux and normal derivatives were studied in depth in Ref.~\cite{Du2013}, yielding a definition of flux $\mathcal{N}$ distinct from the $\mathcal{N}_s$ defined by \eqref{eq:dipierro_neumann_MAIN}; we relate these fluxes in the same appendix.

\vspace*{-8pt}
\subsection{Convergence of truncated fractional operators to untruncated fractional operators}\label{sec:convergence}

Ultimately, we would like to use fractional operators to simulate phenomena such as those mentioned in the introduction. However, discretization of fractional operators is nontrivial and computationally expensive; this is partially due to the fact that one has to deal with infinite domains of integration. The use of truncated operators in place of their untruncated counterpart relaxes the computational challenges by restricting the domains to balls of finite radius in $\mbRn$, while introducing an approximation error. For this reason, it is mandatory to estimate the discrepancy between truncated and untruncated operators. A first step in this direction can be found in Refs. \cite{Defterli2015,DElia2013} where the authors prove the convergence of solutions to the Dirichlet problem for the truncated fractional Laplacian to solutions of the untruncated problem as $\delta\to\infty$. Here, we show that truncated fractional operators converge to their untruncated counterpart in the $L^2$ norm and pointwise. The proofs of the following two theorems, complete with convergence rates, are reported in Appendices \ref{sec:L2convergence} and \ref{sec:point-convergence}.

\vspace*{-4pt}

\begin{theorem}\label{thm:L2convergence}
For $u \in H^s(\mbRn) \cap L^1(\mbRn)$ and $\vb \in {\bf H}^s(\mbRn)$ and $|\vb| \in L^1(\mbRn)$, the truncated fractional divergence and gradient converge in the $L^2$ norm as $\delta\to\infty$ to the fractional divergence and gradient defined in \eqref{eq:frac-Cart}, i.e.
\vskip-10pt %
\begin{equation}\label{eq:norm-conv}
\begin{aligned}
\|\Dsd \vb - \Ds\vb\|_{L^2(\mbRn)}\to 0 & \quad {\rm as}\;\; \delta\to\infty\;\; \forall\;\vb\in [L^2(\mbRn)]^2, \\
\|\Gsd u - \Gs u \|_{L^2(\mbRn)}\to 0    & \quad {\rm as}\;\; \delta\to\infty\;\; \forall\;u\in L^2(\mbRn).
\end{aligned}
\end{equation}
\end{theorem}

\vspace*{-10pt} 

\begin{theorem}\label{thm:point-convergence}
For $u \in H^s(\mbRn)$ and $\vb \in {\bf H}^s(\mbRn)$, the truncated fractional divergence and gradient converge pointwise as $\delta\to\infty$ to the fractional divergence and gradient defined in \eqref{eq:frac-Cart}, i.e.
\vskip -10pt %
\begin{equation}\label{eq:pointwise-conv}
\begin{aligned}
\Dsd \vb(\xb) & \to\Ds\vb(\xb) & \hbox{for almost all $\xb$ as}\;\; \delta\to\infty\;\; \forall\;\xb\in\mbRn, \\
\Gsd u(\xb)   & \to \Gs u(\xb) & \hbox{for almost all $\xb$ as}\;\; \delta\to\infty\;\; \forall\;\xb\in\mbRn.
\end{aligned}
\end{equation}
\end{theorem}

\section{Equivalence of weighted and unweighted \break nonlocal Laplacian operators} \label{sec:unweighted-weighted-equivalence}

\setcounter{section}{4}
\setcounter{equation}{0}\setcounter{theorem}{0}

The theory developed in Section \ref{sec:fractional-nonlocal} follows the framework of weighted nonlocal vector calculus. However, most of the analytical results in the literature, such as well-posedness results \cite{Du2012}, pertain to the unweighted operator $\mathcal{L}$. To leverage these results, a natural question is whether operators of the form $\mathcal{L}_\omega$ can be identified with operators $\mathcal{L}$ of the form \eqref{eq:unw-lapl} with an appropriate kernel $\gamma$.

The ensemble of Theorems \ref{lemma:polar-Cart-equivalence}, \ref{thm:frac-div-grad-lapl}, and \ref{thm:fractional_special_case_nonlocal} provides one example of weighted nonlocal operators, namely $\text{div}^s$ and $\text{grad}^s$, for which $\mathcal{L}_{\omega} = \text{div}^s \text{grad}^s$ can be identified with an unweighted nonlocal Laplacian operator $\mathcal{L}$, namely $-(-\Delta)^{s}$. In this section, we show that such a result holds for general weighted nonlocal operators $\mcD_\omega$ and $\mcG_\omega$. By formally deriving a kernel, we introduce a universal nonlocal Laplacian operator that stems from the composition of weighted divergence and gradient, yet is equivalent to an operator of the form \eqref{eq:unw-lapl}.

We also show that the universal operator corresponds to well-known fractional operators for specific choices of kernels and weights. {Consistently with Theorem \ref{thm:frac-div-grad-lapl}, we show directly that the equivalence kernel for the composition of fractional divergence and gradient evaluates to the fractional Laplacian kernel.} We also derive a representation of the universal operator for tempered weights, suggesting asymptotic equivalence to the tempered fractional Laplacian. These results allow us to consolidate the unified nonlocal vector calculus with the theory introduced by Meerschaert and co-authors in Ref.~\cite{Meerschaert2006}.

\vspace*{-12pt}
\subsection{The equivalence kernel}\label{sec:weighted-unweighted} 

We establish the equivalence of the weighted Laplacian $\mcLw$ and the unweighted Laplacian $\mcL$ for some choice of kernel $\gameq=\gameq(\xb,\yb;\omega,\alphab)$, which we refer to as {\it equivalence kernel}. Recall, once again, that the truncation of the kernel is embedded in $\alphab$ so that both truncated and untruncated cases are accounted for.
\vspace*{-3pt}

\begin{theorem}\label{thm:weighted-unweighted-equivalence}
Let $\mcDw$ and $\mcGw$ be the operators associated with the symmetric weight function $\omega$ and the anti-symmetric function $\alphab$ for $\delta\in(0,\infty]$ and let $u:\Omega\cup\Omega_I^\omega\to\mbR$. For the equivalence kernel $\gameq$ defined by
\vskip -12pt %
\begin{equation}\label{eq:kerndef}
\begin{aligned}
2\gameq(\xb,\yb)
&= \int_\mbRn [\alphab(\xb,\yb)\omega(\xb,\yb)\cdot
   \alphab(\xb,\zb)\omega(\xb,\zb)\\
&  \hspace{1cm}+ \alphab(\zb,\yb)\omega(\zb,\yb)\cdot
   \alphab(\xb,\yb)\omega(\xb,\yb)\\[2mm]
&  \hspace{1cm}+ \alphab(\zb,\yb)
   \omega(\zb,\yb)\cdot\alphab(\xb,\zb)\omega(\xb,\zb)]d\zb,
\end{aligned}
\end{equation}
\vskip -3pt \noindent %
the weighted operator $\mcLw = \mcDw\mcGw$ and unweighted Laplacian operator $\mcL$ with kernel $\gamma_{\textup{eq}}$ are equivalent, i.e.
$
\mcL = \mcLw.
$
\end{theorem}

\proof
By the definition of weighted divergence and gradient and the symmetry of $\omega(\xb,\yb)$, we have
\vskip -12pt %
\begin{align}
   \mcDw\mcGw u(\xb)
&= \int_\mbRn (\mcGw u(\xb) + \mcGw u(\yb))\cdot
   \alphab(\xb,\yb)\omega(\xb,\yb)d\yb \notag\\
&= \int_\mbRn \Bigg[ \int_\mbRn (u(\zb)-u(\xb))
   \alphab(\xb,\zb)\omega(\xb,\zb)d\zb\notag\\
&\qquad\qquad + \int_\mbRn(u(\zb)-u(\yb))
   \alphab(\yb,\zb)\omega(\yb,\zb)d\zb \Bigg] \cdot
   \alphab(\xb,\yb)\omega(\xb,\yb)d\yb \notag\\
&= \int_\mbRn \int_\mbRn (u(\zb)-u(\xb))
   \alphab(\xb,\zb)\omega(\xb,\zb)\cdot
   \alphab(\xb,\yb)\omega(\xb,\yb)d\yb d\zb  
   \\
&\quad+ \int_\mbRn \int_\mbRn (u(\zb)-u(\yb))
   \alphab(\yb,\zb)\omega(\yb,\zb)\cdot
   \alphab(\xb,\yb)\omega(\xb,\yb)d\yb d\zb.
   5\label{eq:eqnB}
\end{align}
\vskip -4pt \noindent
Let the first double integral above be $I$ and the second one be $II$.
We have
\vskip -12pt %
\begin{align*}
I &= \int_\mbRn \int_\mbRn (u(\zb)-u(\xb))
     \alphab(\xb,\zb)\omega(\xb,\zb)\cdot
     \alphab(\xb,\yb)\omega(\xb,\yb)d\yb d\zb \\
  &= \int_\mbRn (u(\zb)-u(\xb))\alphab(\xb,\zb)\omega(\xb,\zb)
     \int_\mbRn \alphab(\xb,\yb)\omega(\xb,\yb)d\yb d\zb.
\end{align*}
\vskip -3pt \noindent%
Letting $\gamma_1(\xb,\zb) = \alphab(\xb,\zb)\omega(\xb,\zb)\int_\mbRn \alphab(\xb,\yb)\omega(\xb,\yb)d\yb$, we have
\vskip -12pt %
\begin{align*}
I &= \int_\mbRn (u(\zb)-u(\xb))\gamma_1(\xb,\zb)d\zb.
\end{align*}
\vskip -6pt \noindent%
Next,
\vskip -14pt
\begin{align*}
II &= \int_\mbRn \int_\mbRn (u(\zb)-u(\yb))
     \alphab(\yb,\zb)\omega(\yb,\zb)\cdot
     \alphab(\xb,\yb)\omega(\xb,\yb)d\yb d\zb \\
  &= \int_\mbRn \int_\mbRn (u(\xb)-u(\yb))
     \alphab(\yb,\zb)\omega(\yb,\zb)\cdot
     \alphab(\xb,\yb)\omega(\xb,\yb)d\yb d\zb \\
  &\qquad+ \int_\mbRn \int_\mbRn (u(\zb)-u(\xb))
     \alphab(\yb,\zb)\omega(\yb,\zb)\cdot
     \alphab(\xb,\yb)\omega(\xb,\yb)d\yb d\zb.
\end{align*}
\vskip -1pt \noindent %
Switching $\yb$ and $\zb$ in the first integral, and employing the anti-symmetry of $\alphab$ and symmetry of $\omega$, we find
\vskip -10pt %
\begin{align*}
II &= \int_\mbRn \int_\mbRn (u(\xb)-u(\zb))
     \alphab(\zb,\yb)\omega(\zb,\yb)\cdot
     \alphab(\xb,\zb)\omega(\xb,\zb)d\zb d\yb\\
  &\qquad+ \int_\mbRn \int_\mbRn (u(\zb)-u(\xb))
     \alphab(\yb,\zb)\omega(\yb,\zb)\cdot
     \alphab(\xb,\yb)\omega(\xb,\yb)d\yb d\zb \\
  &= \int_\mbRn \int_\mbRn (u(\zb)-u(\xb))
     \alphab(\yb,\zb)\omega(\yb,\zb)\cdot
     \alphab(\xb,\zb)\omega(\xb,\zb)d\zb d\yb\\
  &\qquad+ \int_\mbRn \int_\mbRn (u(\zb)-u(\xb))
     \alphab(\yb,\zb)\omega(\yb,\zb)\cdot
     \alphab(\xb,\yb)\omega(\xb,\yb)d\yb d\zb
     \end{align*}
     \begin{align*}
  &= \int_\mbRn \int_\mbRn (u(\zb)-u(\xb))
     \alphab(\yb,\zb)\omega(\yb,\zb)\cdot
     [\alphab(\xb,\zb)\omega(\xb,\zb) + \alphab(\xb,\yb)
     \omega(\xb,\yb)]d\yb d\zb \\
  &= \int_\mbRn \!\!(u(\zb)-u(\xb)) \int_\mbRn
     \alphab(\yb,\zb)\omega(\yb,\zb) \cdot
     [\alphab(\xb,\zb)\omega(\xb,\zb) + \alphab(\xb,\yb)
     \omega(\xb,\yb)]d\yb d\zb.
\end{align*}
Letting $\gamma_2(\xb,\zb)= \int_\mbRn \alphab(\yb,\zb)\omega(\yb,\zb) \cdot [\alphab(\xb,\zb)\omega(\xb,\zb) + \alphab(\xb,\yb)\omega(\xb,\yb)]d\yb$ gives us
\vskip -12pt %
\begin{align*}
II&= \int_\mbRn (u(\zb)-u(\xb))\gamma_2(\xb,\zb)d\zb.
\end{align*}
By combining the above, we have
\begin{align*}
   \mcDw\mcGw u(\xb) &= I+II\\
&= \int_\mbRn (u(\zb)-u(\xb))\gamma_1(\xb,\zb)d\zb
 + \int_\mbRn (u(\zb)-u(\xb))\gamma_2(\xb,\zb)d\zb\\
&= \int_\mbRn (u(\zb)-u(\xb)) (\gamma_1(\xb,\zb)
 + \gamma_2(\xb,\zb))d\zb.
\end{align*}
Then, \eqref{eq:kerndef}, is obtained by setting $2\gameq=\gamma_1+\gamma_2$.
\endproof

\begin{remark}\label{fubini}
In proving Theorem \ref{thm:weighted-unweighted-equivalence}, the order of integration of iterated integrals was repeatedly changed. According to the Fubini-Tonelli theorem, these steps are valid if the iterated integral defining $\mcDw \mcGw u$ converges, or if the integral defining the kernel \eqref{eq:kerndef} as well the integral defined by $\mcLw u$ with that kernel converges. For a given choice of $\alphab$ and $\omega$, it is often simpler to use the latter condition to obtain the appropriate function space in which Theorem \ref{thm:weighted-unweighted-equivalence} is valid.
\end{remark}

\vspace*{-14pt}

\begin{remark}\label{formal_rep_remark}
Theorem \ref{thm:weighted-unweighted-equivalence} only provides a formal equivalence of $\mcL$ and $\mcL_\omega$, in the sense that $\mcL_\omega$ can be represented in the form of equation \eqref{eq:unw-lapl} using $\gamma_{\textup{eq}}$. In fact, as mentioned in Section \ref{sec:unweighted_operators}, the kernel characterizing the unweighted operator $\mcL$ within the theory of nonlocal vector calculus is positive by construction, whereas nothing has been said regarding the sign of the equivalence kernel $\gameq$. Characterization of properties of $\gameq$ is fundamental to analyze well-posedness of nonlocal weighted problems such as \eqref{eq:weight-bound-truncated-prob} and it is discussed later on in this section. However, a full, rigorous characterization of the equivalence kernel and its connection to the well-posedness of the associated problem has not been conducted in the current work, except for a class of singular kernels, see Theorem \ref{eq:beta-equivalence}. As such, our study provides the groundwork for a rigorous unified calculus, which is the subject of current studies.
\end{remark}

\vspace*{-14pt}

\begin{remark}\label{other_eq_kernels}
Instances of the equivalence kernel \eqref{eq:kerndef} have appeared in various contexts in previous works. In the context of peridynamics, Lemma 5 in Ref.~\cite{mengesha2014peridynamic} provided an expression of $\gameq$ for the linear peridynamic Navier equation that is equivalent to \eqref{eq:kerndef}. In that work, the authors conduct a rigorous analysis and provide several results for variational forms associated to truncated weighted operators. For the same peridynamic equation, Ref.~\cite{Alali2015Calculus} also provided a reformulation in terms of a generalized nonlocal calculus that is similar to the one proposed in Ref.~\cite{mengesha2014peridynamic}, and in a related preprint (Ref.~\cite{Alali-preprint}) proposed an instance of $\gameq$ for the linear peridynamic Navier equation.
\end{remark}

\vspace*{-15pt}

\begin{remark}\label{remark:truncation_equivalence}
For $\xb\in\omg$, when we explicitly take into account the truncation of the function $\alphab$, we can rewrite \eqref{eq:kerndef} as
\vskip -13pt %
\begin{align*}
2\gameq(\xb&,\yb) \\
&= \mathds{1}(|\xb-\yb|\leq \delta)
   \rhob(\xb,\yb)\omega(\xb,\yb)\cdot
   \int_{\omg_\delta} \rhob(\xb,\zb)\omega(\xb,\zb)
   \mathds{1}(|\xb-\zb|\leq\delta)\,d\zb\\
&+ \mathds{1}(|\xb-\yb|\leq \delta)
   \rhob(\xb,\yb)\omega(\xb,\yb)\cdot
   \int_{\omg_{2\delta}} \rhob(\zb,\yb)\omega(\zb,\yb)
   \mathds{1}(|\yb-\zb|\leq\delta)\,d\zb\\
&+ \int_{\omg_{2\delta}} \rhob(\zb,\yb)\omega(\zb,\yb)
    \mathds{1}(|\yb-\zb|\leq \delta)\cdot
    \rhob(\xb,\zb)\omega(\xb,\zb)
    \mathds{1}(|\xb-\zb|\leq\delta)\,d\zb
\end{align*}
\vskip -3pt \noindent %
where the last integral is non-zero for $|\xb-\yb|\leq 2 \delta$. Equivalently,
\begin{equation}\label{eq:kerndef-trunc}
\begin{aligned}
2\gameq(\xb,\yb)
&= \mathds{1}(|\xb-\yb|\leq \delta)
   \rhob(\xb,\yb)\omega(\xb,\yb)\cdot
   \int_{B_\delta(\xb)} \rhob(\xb,\zb)\omega(\xb,\zb)\,d\zb\\
&+ \mathds{1}(|\xb-\yb|\leq \delta)
   \rhob(\xb,\yb)\omega(\xb,\yb)\cdot
   \int_{B_\delta(\yb)} \rhob(\zb,\yb)\omega(\zb,\yb)\,d\zb\\
&+ \int_{B_\delta(\xb)\cap B_\delta(\yb)}
    \rhob(\zb,\yb)\omega(\zb,\yb)\cdot
    \rhob(\xb,\zb)\omega(\xb,\zb)\,d\zb.
    \end{aligned}
\end{equation}
Here, the support of the last term in the expression above is $B_{2\delta}(\xb)$; this implies that the interaction domain corresponding to $\mcL$ is a layer of thickness $2\delta$ surrounding $\omg$, as already stated in Section \ref{sec:intro_weighted_operators} and in Refs. \cite{delia2017handbook} and \cite{seleson2016convergence}. Thus, the unweighted operator $\mcL$ associated with the weighted operator $\mcL_\omega$ through the equivalence kernel \eqref{eq:kerndef-trunc} is such that its interaction domain is $\omg_I^\omega=\omg_{2\delta}\setminus\omg$.
\end{remark}

The next corollary shows that when the weight and kernel functions $\omega$ and $\rhob$ are translation invariant, thanks to the anti-symmetry of the integrands, the equivalence kernel reduces to one term. This is particularly useful to show that the equivalence kernel is consistent with well-known radial kernels, such as fractional kernels (see Section \ref{sec:gamma-consistency}).

\vspace*{-3pt}

\begin{corollary}\label{corollary:equiv_kernel_translation}
If the functions $\omega$ and $\rhob$ and  are translation invariant, i.e.,
\vskip -15pt
\begin{equation}
\rho = \rho(\xb - \yb) \quad \text{and} \quad \omega = \omega(\xb - \yb)
\end{equation}
equations \eqref{eq:kerndef} and \eqref{eq:kerndef-trunc} reduce to a single integral. Specifically, \eqref{eq:kerndef} becomes
\vskip -15pt %
\begin{align}\label{eq:kerndef-transl}
2\gameq(\xb,\yb) = \int_\mbRn
\rhob(\zb-\yb)\omega(\zb-\yb)\cdot\rhob(\xb-\zb)\omega(\xb-\zb)\,d\zb,
\end{align}
\vskip -3pt \noindent%
and \eqref{eq:kerndef-trunc} reduces to
\vskip -12pt%
\begin{align}\label{eq:kerndef-trunc-transl}
2\gameq(\xb,\yb) = \int_{B_\delta(\xb)\cap B_\delta(\yb)}
    \rhob(\zb-\yb)\omega(\zb-\yb)\cdot
    \rhob(\xb-\zb)\omega(\xb-\zb)\,d\zb,
\end{align}
\vskip -1pt \noindent%
with support in $B_{2\delta}(\xb)$.
\end{corollary}

\proof
Under these assumptions, the contribution from the first term in \eqref{eq:kerndef} is
\vskip -14pt
\begin{multline}
\alphab(\xb - \yb)\omega(\xb - \yb)\cdot \int_\mbRn \alphab(\xb- \zb)\omega(\xb -\zb) d\zb \\
=
\alphab(\xb - \yb)\omega(\xb - \yb)\cdot \int_\mbRn \alphab(\zb')\omega(\zb') d\zb.
\end{multline}
\vskip -3pt \noindent
The integral on the right-hand side is zero due to antisymmetry of $\alphab(\zb')\omega(\zb')$. Similarly, one can show the contribution from the second term in \eqref{eq:kerndef} is zero. The contribution from the third term in \eqref{eq:kerndef} is then identical to \eqref{eq:kerndef-transl}, and \eqref{eq:kerndef-trunc-transl} can be obtained just as in Remark \ref{remark:truncation_equivalence}.
\endproof

\vspace*{-10pt}
\subsection{Properties of the equivalence kernel}
Expression \eqref{eq:kerndef} establishes an important relation between the weight and kernel functions associated with $\mcLw$ and the kernel function $\gamma_\textup{eq}$ yielding a representation of $\mcLw$ in the form \eqref{eq:unw-lapl} of an operator $\mcL$. However, as pointed out in Remark \ref{formal_rep_remark}, this relation is not sufficient to transfer well-posedness results from the unweighted to the weighted case.

The following two lemmas prove fundamental properties of the equivalence kernel and provide conditions on $\omega$ and $\alphab$ that will guarantee coercivity of the weighted operator via the equivalence with $\mcL$, as explained in Section \ref{sec:weak-equivalence}. 

\vspace*{-4pt}

\begin{lemma}
The equivalence kernel $\gameq$ in \eqref{eq:kerndef} is symmetric for all $(\xb,\yb)\in\mbRn$ and for $\delta\in(0,\infty]$.
\end{lemma}
\proof
 This follows from the anti-symmetry of $\alphab$ and the symmetry of $\omega$. We let $\etab(\xb,\yb)$ be the anti-symmetric function defined as $\etab(\xb,\yb)=\alphab(\xb,\yb)\omega(\xb,\yb)$; then, we rewrite \eqref{eq:kerndef} as
\vskip -11pt%
\begin{align*}
2\gameq(\xb,\yb) &=  \int_\mbRn [\etab(\xb,\yb)\cdot \etab(\xb,\zb) + \etab(\zb,\yb)\cdot\etab(\xb,\yb) + \etab(\zb,\yb)\cdot\etab(\xb,\zb)]d\zb
\end{align*}
The anti-symmetry of $\etab$ implies that
\vskip-11pt
\begin{align*}
2\gameq(\xb,\yb)&= \int_\mbRn [\etab(\yb,\xb)\cdot \etab(\zb,\xb) + \etab(\yb,\zb)\cdot\etab(\yb,\xb) + \etab(\yb,\zb)\cdot\etab(\zb,\xb)]d\zb
\end{align*}
Since the dot product is commutative, we switch the orders of each $\etab$ pair:
\vskip -10pt%
\begin{align*}
2\gameq(\xb,\yb)&= \int_\mbRn [\etab(\zb,\xb)\cdot\etab(\yb,\xb) + \etab(\yb,\xb)\cdot\etab(\yb,\zb) + \etab(\zb,\xb)\cdot\etab(\yb,\zb)]d\zb
\end{align*}
then, switching the first two terms, we have
\begin{align*}
2\gameq(\xb,\yb) &= \!\int\limits_\mbRn [\etab(\yb,\xb)\cdot\etab(\yb,\zb) + \etab(\zb,\xb)\cdot\etab(\yb,\xb) + \etab(\zb,\xb)\cdot\etab(\yb,\zb)]d\zb\\ &=2\gameq(\yb,\xb).
\end{align*}\vspace*{-7pt}
\endproof 

The following theorem shows that when $\rhob$ is the unit vector in the direction $\xb-\yb$ and $\omega$ is a power of the modulus of $\xb-\yb$, then we can write the equivalence kernel in a closed form and show its positivity. Kernels of this type are often encountered in the fractional setting.

\vspace*{-3pt}

\begin{theorem}\label{eq:beta-equivalence}
Let $\beta>n/2$ and $\omega$ and $\rhob$ be defined as
\vskip -12pt%
\begin{equation}\label{eq:omega-alpha-wellposed}
\begin{aligned}
\omega(\xb,\yb) &= |\yb-\xb|\phi(|\yb-\xb|)
\;\; \text{with} \;\;
\phi(|\yb-\xb|)=|\yb-\xb|^{-\beta}\\
\rhob(\xb,\yb) &= \frac{\yb-\xb}{|\yb-\xb|}
\;\; \mbox{with} \;\;
\delta=\infty
\end{aligned}
\end{equation}
\vskip -2pt \noindent%
then, the equivalence kernel $\gameq(\xb,\yb;\omega,\alphab)$ defined in \eqref{eq:kerndef} is such that
\vskip -9pt%
\begin{equation}\label{eq:gamma-wellposed}
\gameq(\xb,\yb)=\gameq(|\yb-\xb|)=
\overline\gamma |\yb-\xb|^{n+2(1-\beta)},
\end{equation}
\vskip -3pt \noindent%
where the constant $\overline\gamma$ is given by
\vskip -9pt
\begin{equation}
2\overline\gamma = \int_\mbRn \dfrac{({\bf e}-\zb)\cdot\zb}
{|{\bf e}-\zb|^\beta|\zb|^\beta}\,d\zb,
\end{equation}
\vskip -3pt \noindent
for any unit vector $\bf e$.
\end{theorem}

\proof
Since $\rhob$ and $\omega$ are translation invariant, according to Corollary \eqref{corollary:equiv_kernel_translation}, we rewrite the equivalence kernel as
\vskip-12pt%
\begin{align*}
2\gameq(\xb,\yb) =
\int_{\mbRn}\frac{\yb-\zb}{|\yb-\zb|^\beta}\cdot
\frac{\zb-\xb}{|\zb-\xb|^\beta} d\zb.
\end{align*}
\vskip -3pt \noindent%
We evaluate this integral indirectly. Let $\zb' = \zb-\xb$. Then $\zb = \zb' + \xb$, $d\zb = d\zb'$ and $\yb - \zb = \yb - \zb' - \xb = \yb - \xb - \zb'$. Thus,
\vskip -12pt%
\begin{align*}
2\gameq(\xb,\yb) &= \int_{\mbRn}\frac{\yb-\xb-\zb'}
{|\yb-\xb-\zb'|^\beta}\cdot
\frac{\zb'}{|\zb'|^\beta} d\zb'\\
&=\int_{\mbRn}\frac{\yb-\xb-\zb}{|\yb-\xb-\zb|^\beta}
\cdot\frac{\zb}{|\zb|^\beta} d\zb.
\end{align*}
\vskip -3pt \noindent%
Thus, $\gameq(\xb,\yb)$ depends only on $\xb-\yb$. Next, we show that $\gameq$ is rotationally invariant. Consider a rotation $\mcR$; we have
\vskip-12pt %
\begin{align*}
2\gameq\left(\mathcal{R}(\xb - \yb)\right) &=
\int_{\mbRn}
\frac{\mathcal{R}(\yb-\xb)-\zb}{|\mathcal{R}(\yb-\xb)-\zb|^\beta}
\cdot \frac{\zb}{|\zb|^\beta} d\zb.
\end{align*}
\vskip -3pt \noindent%
Let $\zb = \mathcal{R} \zb'$. Then $d\zb = d\zb'$, and
\begin{align*}
2\gameq\left(\mathcal{R}(\xb-\yb)\right)
&= \int_{\mbRn}\frac{\mathcal{R}(\yb-\xb)-\mathcal{R}\zb'}
{|\mathcal{R}(\yb-\xb)-\mathcal{R}\zb'|^\beta}
\cdot\frac{\mathcal{R}\zb'}
{|\mathcal{R}\zb'|^\beta} d\zb'
\\
&=
\int_{\mbRn}\frac{\mathcal{R}(\yb-\xb)-\mathcal{R}\zb}
{|\mathcal{R}(\yb-\xb)-\mathcal{R}\zb|^\beta}
\cdot\frac{\mathcal{R}\zb}{|\mathcal{R}\zb|^\beta} d\zb
\\
&=
\int_{\mbRn}\frac{\mathcal{R}\left((\yb-\xb)-\zb\right)}
{|\mathcal{R}\left((\yb-\xb)-\zb\right)|^\beta}\cdot
\frac{\mathcal{R}\zb}{|\mathcal{R}\zb|^\beta} d\zb
\\
&=
\int_{\mbRn}
\frac{(\yb-\xb)-\zb}{|(\yb-\xb)-\zb|^\beta}
\cdot\frac{\zb}{|\zb|^\beta}\ d\zb\\
&=2\gameq(\xb-\yb).
\end{align*}
Therefore, $\gameq(\xb,\yb)$ depends only on $|\xb-\yb|$. Now we let $c > 0$ and consider
\vskip -10pt%
\begin{align*}
2\gameq(c|\xb-\yb|)
&= \int_{\mbRn} \frac{c(\yb-\xb)-\zb}{|c(\yb-\xb)-\zb|^\beta}
\cdot \frac{\zb}{|\zb|^\beta} d\zb.
\end{align*}
Let $\zb = c \zb'$. Then $d\zb = c^n d\zb'$, and
\begin{align*}
2\gameq(c|\xb-\yb|)
&= \int_{\mbRn}\frac{c(\yb-\xb)-c \zb'}
{|c(\yb-\xb)-c \zb'|^\beta}
\cdot\frac{c \zb'}{|c \zb'|^\beta}
c^n d\zb' \\
&=
\int_{\mbRn}
\frac{c(\yb-\xb)-c \zb}
{|c(\yb-\xb)-c \zb|^\beta}
\cdot\frac{c \zb}{|c \zb|^\beta}c^ndz \\
&=
c^{n+2(1-\beta)}
\int_{\mbRn}\frac{(\yb-\xb)-\zb}{|(\yb-\xb)-\zb|^\beta}
\cdot\frac{\zb}{|\zb|^\beta}dz \\
&=
c^{n+2(1-\beta)} 2\gameq(|\xb-\yb|).
\end{align*}
\vskip -3pt \noindent%
Therefore,
\vskip -12pt
\begin{align*}
2\gameq(\xb,\yb) &= |\xb-\yb|^{n+2(1-\beta)} 2\gameq({\bf e}),
\end{align*}
where $\bf e$ is any unit vector and $\gameq({\bf e})=\overline\gamma$ is a constant independent of the choice of unit vector $\bf e$:
\vskip -12pt%
\begin{equation*}
2\gameq({\bf e})=\int_{\mbRn}\frac{{\bf e}-\zb}{|{\bf e}-\zb|^\beta}
\cdot \frac{\zb}{|\zb|^\beta} d\zb.
\end{equation*}
\vspace*{-4pt}
\endproof

\vspace*{-10pt}
\subsection{Consistency of the equivalence for fractional-type operators}\label{sec:gamma-consistency} 

In this section we show that for specific choices of $\omega$ and $\rhob$ the equivalence kernel $\gameq$ corresponds to well-known fractional kernels. The following result is a corollary to Theorem \ref{eq:beta-equivalence} and proves the consistency of the equivalence kernel for the fractional Laplacian operator.

\vspace*{-3pt}

\begin{corollary}[The fractional Laplacian kernel]
\label{thm:frac-equivalence-kernel}
For $\omega$ and $\rhob$ defined as
\vskip -14pt%
\begin{align*}
\omega(\xb,\yb) &= |\yb-\xb|\phi(|\yb-\xb|)
\;\; \text{with} \;\;
\phi(|\yb-\xb|)=\frac{1}{|\yb-\xb|^{n+1+s}}\\
\rhob(\xb,\yb) &= \frac{\yb-\xb}{|\yb-\xb|}
\;\; \text{with} \;\; \delta=\infty
\end{align*}
\vskip -3pt \noindent %
for which \,
$\alphab(\xb,\yb)\omega(\xb,\yb) = (\yb-\xb)|\yb-\xb|^{-(n+s+1)}$, we have
\vskip -8pt%
\begin{equation}\label{eq:frac-kernel-equivalence}
2\gameq(\xb,\yb)=-\frac{C_{n,s}D_{n,s}}{G_s^2}
\dfrac{1}{|\xb-\yb|^{n+2s}},
\end{equation}
where $G_s=s/\Gamma(1-s)$ and $D_{n,s}$ is the negative constant defined as in Lemma \ref{lem:integral-equality-constant}. Then, Theorem \ref{thm:weighted-unweighted-equivalence} is valid for $u \in H^s(\mathbb{R}^n)$, and
\vskip -8pt
\begin{equation}
\mcL= \mcLw = {\frac{D_{n,s}}{G_s^2}}(-\Delta)^s.
\end{equation}
\end{corollary}

\proof
The $\phi$ defined above are the same form as in Theorem \ref{eq:beta-equivalence} with
\vskip -15pt
\begin{equation}
\beta=n+1+s>n/2,
\end{equation}
\vskip -1pt \noindent
so the assumptions of Theorem \ref{eq:beta-equivalence} are satisfied and $\gamma_{\text{eq}}$ converges to the stated expression. Therefore, for constant $\overline\gamma$,
\vskip -9pt
\begin{equation}
\mcLw u (\xb)=
2\overline\gamma
\int_{\mathbb{R}^d}
\frac{u(\yb) - u(\xb)}
{|\yb-\xb|^{n+2s}} d\yb,
\end{equation}
which converges, in a principal value sense, for $u \in H^s(\mathbb{R}^n)$.
On the other hand, from Theorems \ref{lemma:polar-Cart-equivalence}, \ref{thm:frac-div-grad-lapl}, and \ref{thm:fractional_special_case_nonlocal}, we have
\vskip -10pt
\begin{equation}
\mcLw u (\xb)=
\frac{D_{n,s}}{G_s^2}
{
(-\Delta)^s u(x)
=
-\frac{C_{n,s} D_{n,s}}{G_s^2}
\int_{\mathbb{R}^d}
\frac{u(\yb) - u(\xb)}
{|\yb-\xb|^{n+2s}} d\yb
}.
\end{equation}
\vskip -2pt \noindent
Therefore $2\overline\gamma = {-C_{n,s} D_{n,s}}/{G_s^2}.$
\endproof

The next theorem derives a representation of $\gamma_{\textup{eq}}$ for tempered fractional weights, corresponding to the tempered fractional gradient and divergence \eqref{eq:tempered-fractional}, which allows for comparison to the kernel of the tempered fractional Laplacian \eqref{eq:tempered_Laplacian_rn}.
\vspace*{-4pt}

\begin{theorem}[The tempered fractional Laplacian kernel] \label{thm:tempered-equivalence-kernel}
For $\omega$ and $\rhob$ defined as
\vskip -14pt%
\begin{align*}
\omega(\xb,\yb) &= |\yb-\xb|\phi(|\yb-\xb|)
\;\; \text{with} \;\;
\phi(|\yb-\xb|)=
\frac{e^{-\lambda|\yb-\xb|}}{|\yb-\xb|^{n+1+s}}
\\
\rhob(\xb,\yb) &= \frac{\yb-\xb}{|\yb-\xb|},
\;\; \text{with} \;\;
\delta=\infty
\end{align*}
\vskip -3pt \noindent %
for which
$\alphab(\xb,\yb)\omega(\xb,\yb) =
{(\yb-\xb)}{|\yb-\xb|^{-(n+1+s)}}e^{-\lambda|\yb-\xb|}
$, we have
\vskip -11pt%
\begin{align}
\label{eq:tempered-kernel-equivalence}
2\gameq(\xb,\yb) &
= \frac{F(n,s,\lambda, |\xb - \yb|)}{|\xb-\yb|^{n+2s}},
\end{align}
\vskip -4pt \noindent
where
\vskip -13pt
\begin{equation}
\label{tempered_equivalence_factor}
F(n, s, \lambda, |\xb - \yb|) =
\int_{\mbRn}
\frac{{\bf e}-\zb}
{|{\bf e}-\zb|^{n+s+1}}
\cdot
\frac{\zb}
{|\zb|^{n+s+1}}
e^{-\lambda|\xb-\yb|\left(|{\bf e - \zb}|+|\zb|\right)}
d\zb.
\end{equation}
\end{theorem}

The proof of this theorem is similar to that of Theorem \ref{thm:frac-equivalence-kernel}, and is given in Appendix \ref{sec:proof_tempered_equivalence_kernel}.

\smallskip 

The above theorem extracts the dependence on $|\xb - \yb|$ from $\gamma_{\textup{eq}}$ in \eqref{eq:tempered-kernel-equivalence}, just as for the standard fractional Laplacian in Corollary \ref{thm:frac-equivalence-kernel}. However, nothing has been proven about the factor $F(n, s, \lambda, |\xb - \yb|)$ given by \eqref{tempered_equivalence_factor}. The kernel \eqref{eq:tempered-kernel-equivalence} should be compared to the kernel of the tempered fractional Laplacian \eqref{eq:tempered_Laplacian_rn}. While we do not expect the factor \eqref{tempered_equivalence_factor} to reduce to $e^{-\lambda|\xb - \yb|}$, the form of \eqref{tempered_equivalence_factor} suggests that
\vskip -11pt
\begin{equation}
F(n, s, \lambda, |\xb - \yb|) \sim e^{-\lambda|\xb - \yb|}
\end{equation}
\vskip -2pt \noindent
in an asymptotic sense. This would imply that the use of tempered fractional divergence and gradient \eqref{eq:tempered-fractional} in the unified vector calculus framework yields an operator $\mathcal{L}_\omega$ consistent in an asymptotic sense with $-(-\Delta)_\lambda^s$. Moreover, establishing non-negativity of $F(n, s, \lambda, |\xb - \yb|)$ as a function of $|\xb - \yb|$ would imply that the resulting exterior-value problem is well-posed, by the arguments of Section \ref{sec:unified-variational}. Thus, we state the following conjecture, which is supported by numerical evidence\cite{Olson2020CSRI}; we leave the full proof for future study.

\vspace*{-3pt}

\begin{conjecture}
{\sl We conjecture that $F$ satisfies the following bounds
\vskip -10pt %
\begin{equation}\label{tempered_equivalence_bounds}
\underline{B}_{n,s} e^{-\lambda|\xb - \yb|} \le F(n, s, \lambda, |\xb - \yb|) \le \overline{B}_{n,s} e^{-\lambda|\xb - \yb|}
\end{equation}
\vskip -2pt \noindent %
for positive constants $\underline{B}_{n,s}$ and $\overline{B}_{n,s}$.
}
\end{conjecture} 

\section{A unified variational setting}  
\label{sec:unified-variational}

\setcounter{section}{5}
\setcounter{equation}{0}\setcounter{theorem}{0}

The results presented in the previous section establish the equivalence of weighted and unweighted operators for a specific choice of kernel $\gameq(\xb,\yb;\omega,$ $\alphab)$. This enables the extension of the well-established analysis of the unweighted Dirichlet problem to the weighted case and, hence, to various fractional cases. In this section we first establish a general Green first identity for weighted operators and then use it to derive weighted weak forms and study their well-posedness. We recall that Refs. \cite{Du2018Dirichlet,Du2020SPH,Lee2020} investigated the well-posedness of problems associated to weighted operators for truncated kernels. Here, we extend the discussion to infinite-range kernels and provide a well-posedness result for the class of kernels introduced in Theorem \ref{eq:beta-equivalence}.

\vspace*{-12pt}
\subsection{A weighted nonlocal Green identity} \label{sec:weighted-Green}

The following theorem proves a weighted Green identity; this result enables the derivation of weighted variational forms.
\vspace*{-3pt}

\begin{theorem}[Nonlocal weighted Green's identity] \label{thm:weighted-Green}
Let $u,v:\Omega\cup\Omega_I^\omega\to\mbR$ and let $\delta\in(0,\infty)$, the weighted nonlocal operators satisfy a (weighted nonlocal) Green's identity, i.e.
\vskip -12pt%
\begin{equation}\label{eq:weighted-Green}
  \int_{\Omega}-\mcLw u(\xb) v(\xb) \ d\xb
= \int_\mbRn \mcGw u(\xb) \!\cdot\! \mcGw u(\xb) d\xb
+ \int_{\mbRn \setminus \Omega} \mcDw\mcGw u(\xb) v(\xb) \ d\xb.
\end{equation}
\end{theorem}

\proof
By definition of weighted divergence \eqref{eq:w-div} and gradient \eqref{eq:w-grad} we have
\vskip -12pt%
\begin{align*}
\begin{split}
\int_{\mbRn} &\mathcal{D}_\omega\mathcal{G}_\omega  u(\xb) v(\xb) \ d\xb \\
&=
\int_{\mbRn}\mathcal{D}\left(\omega(\xb,\yb)\mathcal{G}_\omega u(\xb)\right)
v(\xb) \ d\xb\\
&=
\int_{\mbRn}\int_{\mbRn}\left[\omega(\xb,\yb) \mathcal{G}_\omega u(\xb)
+\omega(\yb,\xb)\mathcal{G}_\omega u(\yb) \right] \cdot
\alphab(\xb,\yb) v(\xb) \ d\yb d\xb\\
&=\begin{multlined}[t]
\int_{\mbRn}\int_{\mbRn}\omega(\xb,\yb)
\left[\int_{\mbRn}\mathcal{G} u(\xb,\zb) \omega(\xb,\zb)d\zb
+\int_{\mbRn}\mathcal{G} u(\yb,\zb) \omega(\yb,\zb)d\zb\right]
\\
\times\, \alphab(\xb,\yb)v(\xb) \ d\yb d\xb.
\end{multlined}\\
\end{split}
\end{align*}
\vskip -3pt \noindent
By separating into two integrals, we have
$$
\int_{\mbRn}
\mathcal{D}_\omega
\mathcal{G}_\omega
u(\xb)  v(\xb) \ d\xb 
\hspace*{7cm}
$$

\begin{align*}
\begin{split}
&=
\begin{multlined}[t]
\int_{\mbRn}
\int_{\mbRn}
\omega(\xb,\yb)
\int_{\mbRn}
\mathcal{G} u(\xb,\zb) \omega(\xb,\zb)
d\zb
\cdot
\alphab(\xb,\yb)
v(\xb)  d\yb d\xb
\\
+\int_{\mbRn}
\int_{\mbRn}
\omega(\xb,\yb)
\int_{\mbRn}
\mathcal{G} u(\yb,\zb) \omega(\yb,\zb)
d\zb
\cdot
\alphab(\xb,\yb)
v(\xb)  d\yb d\xb.
\end{multlined}
\end{split}
\end{align*}
\begin{align*}
\begin{split}
&=
\begin{multlined}[t]
\int_{\mbRn}
\int_{\mbRn}
\int_{\mbRn}
\omega(\xb,\yb)
\omega(\xb,\zb)
\mathcal{G} u(\xb,\zb)
\cdot
\alphab(\xb,\yb)
v(\xb)
d\zb d\yb d\xb \\
+
\int_{\mbRn}
\int_{\mbRn}
\int_{\mbRn}
\omega(\xb,\yb)
\omega(\yb,\zb)
\mathcal{G} u(\yb,\zb)
\cdot
\alphab(\xb,\yb)
v(\xb)
d\zb d\yb d\xb.
\end{multlined}
\end{split}
\end{align*}
The change of variables $(\xb,\yb,\zb) \rightarrow (\yb,\zb,\xb)$ in the second term gives
\vskip -11pt%
\begin{align*}
\int_{\mbRn}
&\mathcal{D}_\omega
\mathcal{G}_\omega
u(\xb)  v(\xb) \ d\xb \\
&=
\begin{multlined}[t]
\int_{\mbRn}
\int_{\mbRn}
\int_{\mbRn}
\omega(\xb,\yb)
\omega(\xb,\zb)
\mathcal{G} u(\xb,\zb)
\cdot
\alphab(\xb,\yb)
v(\xb)
d\zb d\yb d\xb
\\+
\int_{\mbRn}
\int_{\mbRn}
\int_{\mbRn}
\omega(\yb,\zb)
\omega(\zb,\xb)
\mathcal{G} u(\zb,\xb)
\cdot
\alphab(\yb,\zb)
v(\yb)
d\xb d\zb d\yb
\end{multlined}
\\
&=
\begin{multlined}[t]
\int_{\mbRn}
\int_{\mbRn}
\omega(\xb,\yb)
\mathcal{G}_\omega u(\xb)
\cdot
\alphab(\xb,\yb)
v(\xb)
d\yb d\xb
\\+
\int_{\mbRn}
\int_{\mbRn}
\omega(\yb,\zb)
\mathcal{G}_\omega u(\zb)
\cdot
\alphab(\yb,\zb)
v(\yb)
 d\zb d\yb
\end{multlined}
\end{align*}
\begin{align*}
&=
\begin{multlined}[t]
\int_{\mbRn}
\int_{\mbRn}
\omega(\xb,\yb)
\mathcal{G}_\omega u(\xb)
\cdot
\alphab(\xb,\yb)
v(\xb)
d\yb d\xb
\\+
\int_{\mbRn}
\int_{\mbRn}
\omega(\yb,\xb)
\mathcal{G}_\omega u(\xb)
\cdot
\alphab(\yb,\xb)
v(\yb)
 d\xb d\yb
\end{multlined}
\\
&=
\int_{\mbRn}
\int_{\mbRn}
\omega(\xb,\yb)
\mathcal{G}_\omega u(\xb)
\cdot
\alphab(\xb,\yb)
\left[
v(\xb)
-
v(\yb)
\right]
 d\xb d\yb
\\
&= -
\int_{\mbRn}
\int_{\mbRn}
\omega(\xb,\yb)
\mathcal{G}_\omega u(\xb)
\cdot
\mathcal{G} v(\xb,\yb)
 d\xb d\yb
\\
&= -
\int_{\mbRn}
\mathcal{G}_\omega u(\xb)
\cdot
\mathcal{G}_\omega v(\xb)
d\xb.
\end{align*}
\vspace*{-2pt}
\endproof 

For $\delta<\infty$  and $\alphab$ with support in $\Omega_{2\delta} \times \Omega_{2\delta}$, both $\mathcal{G}_\omega u(\xb)$ and $\mathcal{D}_\omega \mathcal{G}_\omega u(\xb)$ are zero for $\xb \in \mathbb{R}^n \setminus \Omega_{2\delta}$, see \eqref{eq:w-grad} and \eqref{eq:w-lapl}. Therefore, under these assumptions Theorem \ref{thm:weighted-Green} yields as a special case
\begin{equation*}
  \int_\omg - \mcLw u(\xb) v(\xb) \, d\xb
= \int_{\omg\cup\omg_I^\omega} \mcGw u(\xb) \!\cdot\! \mcGw u(\xb) d\xb
+ \int_{\omg_I^\omega} \mcDw\mcGw u(\xb) v(\xb) \, d\xb
\end{equation*}
This result, i.e. a Green identity for truncated weighted operators, was stated without proof in Ref.~\cite{du2013analysis}. The same result was stated and proved in Ref.~\cite{mengesha2016characterization} with slightly different conventions and definitions for the weighted operators in Section \ref{sec:intro_weighted_operators}, which amount to the same assumption that $\alphab$ has support in $\Omega_{2\delta} \times \Omega_{2\delta}$. The identification in Section \ref{sec:fractional-nonlocal-equivalence} of fractional gradient and fractional divergence as special cases of weighted nonlocal operators, together with the identification of their composition with an appropriate nonlocal diffusion operator in Theorem \ref{thm:weighted-unweighted-equivalence}, means that Theorem \ref{thm:weighted-Green} provides a novel Green identity for such operators, in addition to their tempered variants.

\vspace*{-6pt}
\subsection{Equivalence of weighted and unweighted variational forms} 
\label{sec:weak-equivalence}

Theorem \ref{thm:weighted-Green} allows us to rewrite \eqref{eq:weight-truncated-weak} in terms of a weighted energy. Green's identity in \eqref{eq:weighted-Green} yields
\vskip -13pt%
\begin{equation}\label{eq:weight-truncated-weak-green}
\int_\Omega (-\mcLw u -f)\, v \,d\xb =
\int_{\omg\cup\omg_I^\omega} \mcGw u \,\mcGw v \,d\xb
- \int_\Omega f\, v \,d\xb = 0,
\end{equation}
\vskip -2pt \noindent %
or, equivalently, $\mcA_\omega(u,v) = \mcF(v),$
where $\mcF(\cdot)$ is defined as in \eqref{eq:A-F-V} and
\vskip -10pt %
\begin{displaymath}
\mcA_\omega(u,v)  =\int_{\omg\cup\omg_I^\omega} \mcGw u \,\mcGw v \,d\xb.
\end{displaymath}
\vskip -3pt \noindent%
Note that $\mcA_\omega(v,v)=\vertiii{v}^2_\omega$, the weighted energy introduced in \eqref{eq:weighted-energy}. However, since we did not prove that the weighted energy defines a norm in the functional space $V(\Omega\cup\Omega_I^\omega)$, this fact is not enough to prove the well-posedness of the weak weighted problem \eqref{eq:weight-truncated-weak-green} for any choice of weight $\omega$ and kernel $\alphab$. To this end, we first prove the equivalence of weighted and unweighted bilinear forms for $\omega$ and $\alphab$ such that \eqref{eq:kerndef} is satisfied.
\vspace*{-3pt}

\begin{theorem}[Variational equivalence]\label{thm:var-equivalence}
For $u,v:\Omega\cup\Omega_I^\omega\to\mbR$, $\delta\in(0,\infty)$ and $\gameq(\xb,\yb;\omega,\alphab)$ defined as in \eqref{eq:kerndef}, the variational forms associated with truncated weighted and truncated unweighted nonlocal operators are equivalent. That is,
\vskip -12pt%
\begin{equation}\label{eq:weak-equivalence}
 \mcA(u,v)
=\int_{\omg\cup\omg_I^\omega}\int_{\omg\cup\omg_I^\omega}
 \mcG u \,\mcG v \,d\yb\,d\xb
=\int_{\omg\cup\omg_I^\omega} \mcGw u \,\mcGw v \,d\xb
=\mcA_\omega(u,v),
\end{equation}
for all $v=0$ in $\omg_I^\omega$, where $\omg_I^\omega$ is the interaction domain associated with the weighted Laplacian operator $\mcLw$.
\end{theorem}

\proof
The proof follows from Theorems \ref{thm:weighted-unweighted-equivalence} and \ref{thm:weighted-Green}. We have
\begin{displaymath}
\begin{aligned}
&     \int_{\omg\cup\omg_I^\omega}\int_{\omg\cup\omg_I^\omega}
      \mcG u\, \mcG v \,d\yb\,d\xb
&     \quad\hbox{(unweighted Green's identity \eqref{eq:unweighted-Green})} \\[2mm]
& = - \int_\omg \mcD\mcG u v \,d\xb
&     \quad\hbox{(Equivalence Theorem \ref{thm:weighted-unweighted-equivalence})} \\[2mm]
& = - \int_\omg \mcDw\mcGw u v \,d\xb
&     \quad\hbox{(weighted Green's identity \eqref{eq:weighted-Green})}\\[2mm]
& =   \int_{\omg\cup\omg_I^\omega} \mcGw u \, \mcGw v \,d\xb. &
\end{aligned}
\end{displaymath}
\vspace*{-5pt}
\endproof 

Once again, note that the case $\delta=\infty$ is also included. An immediate consequence of this theorem is the formal equivalence of the weighted and unweighted energies; in fact,
\begin{equation}\label{eq:equivalence_of_norms}
\vertiii{v}^2 = \mcA(v,v) = \mcA_\omega(v,v)= \vertiii{v}^2_\omega.
\end{equation}

We are now ready to state the weak form of problem \eqref{eq:weight-bound-truncated-prob} and prove its well posedness. For $f\in V'(\omg\cup\omg_I^\omega)$, find $u\in V(\omg\cup\omg_I^\omega)$ such that $u=0$ in $\omg_I^\omega$ and
\vskip -10pt %
\begin{equation}\label{eq:weighted-weak-AF}
\mcA_\omega(u,v) = \mcF(v), \;\;\forall\, v\in V_c(\omg\cup\omg_I^\omega).
\end{equation}
The following theorem proves the well posedness of \eqref{eq:weighted-weak-AF} for a weighted Laplacian whose equivalence kernel is either positive or negative. Note that this theorem gives closure to Open Problem 1.9 in Ref. \cite{Shieh2015} for the case of a constant, scalar, diffusion parameter. For the extension to anisotropic tensor-valued diffusivity, see Ref. \cite{d2021analysis}.

\vspace*{-3pt}

\begin{theorem}\label{thm:weighted-well-posedness}
If the equivalence kernel \eqref{eq:kerndef} for the operator $\mcL_\omega$ corresponding to the the bilinear form $\mcA_\omega$ satisfies $\gamma_{\textup{eq}}(\xb,\yb) \geq 0$ or $\gamma_{\textup{eq}}(\xb,\yb) \leq 0$ for all $\xb, \yb$, then problem \eqref{eq:weighted-weak-AF} is well-posed in $V_c(\omg\cup\omg_I^\omega)$.
\end{theorem}
\proof
Recall that the equivalence kernel $\gameq(\xb,\yb;\omega,\alphab)$ in \eqref{eq:kerndef} is symmetric. We first consider the nonnegative case. If $\gameq\geq 0$ for all $\xb, \yb$, then the associated unweighted bilinear form $\mcA$ in \eqref{eq:A-F-V} is coercive. Due to Theorem \ref{thm:var-equivalence}, the weighted bilinear form $\mcA_\omega$ is coercive and continuous in $V(\omg\cup\omg_I^\omega)$. This fact and the continuity of $\mcF$ guarantees the well-posedness of \eqref{eq:weighted-weak-AF} by the Lax-Milgram theorem. If the kernel is nonpositive everywhere, by simply multiplying the strong form \eqref{eq:weight-bound-truncated-prob} of the nonlocal equation by $-1$, the well-posedness result for nonnegative kernels applies. \break \hspace*{1cm} 
\endproof
\vspace*{-1pt}

A direct consequence of Theorem \ref{thm:weighted-well-posedness} is that, for $\omega$ and $\alphab$, defined as in \eqref{eq:omega-alpha-wellposed} problem \eqref{eq:weighted-weak-AF} is well-posed. The proof simply follows from the fact that, with these choices, $\gameq(\xb,\yb;\omega,\alphab)$ in \eqref{eq:kerndef} is a constant times a positive function; as such, it satisfies assumptions of Theorem \ref{thm:weighted-well-posedness}.

\vspace*{-3pt}

\begin{corollary}
If the equivalence kernel \eqref{eq:kerndef} for the operator $\mcL_\omega$ corresponding to the the bilinear form $\mcA_\omega$ is such that $\omega$ and $\alphab$ are defined as in \eqref{eq:omega-alpha-wellposed}, then problem \eqref{eq:weighted-weak-AF} is well-posed in $V_c(\omg\cup\omg_I^\omega)$.
\end{corollary}

\vspace*{-10pt}

\begin{remark}
Theorem \ref{thm:weighted-well-posedness} provides a condition for proving well-posedness of \eqref{eq:weighted-weak-AF} by analyzing the equivalence kernel.
The above corollary demonstrates that this yields well-posedness for certain power-law kernels, including those corresponding to fractional-order operators. While this strategy may be utilized for other types of kernels, problem \eqref{eq:weighted-weak-AF} may fail to be well-posed if the equivalence kernel is not sufficiently singular \cite{Du2018Dirichlet}.
This can be resolved in certain so-called \emph{correspondence models} of peridynamics by the addition of terms in the equation to obtain a well-posed problem; see Refs. \cite{mengesha2014peridynamic,silling2017stability}.
On the other hand, sufficient conditions on the singularity of the kernel that guarantee well-posedness have been obtained in Refs. \cite{Du2018Dirichlet,Du2020SPH,Lee2020}. A further study is still needed to fully characterize classes of kernels whose associated problems are well-posed.
\end{remark}

\vspace*{-7pt}
\section{Other instances of fractional vector calculus}\label{sec:Tarasov}

\setcounter{section}{6}
\setcounter{equation}{0}\setcounter{theorem}{0}

During the last decade, several distinct theories of fractional vector calculus have appeared in the literature that are based on combinations of one-dimensional fractional derivatives with respect to individual coordinates.
Ref. \cite{Tarasov2008} provides a review of some of those frameworks, while also introducing and developing such a theory based on one-dimensional Caputo derivatives. Other examples of such theories are that of Adda \cite{Adda1997}, based on one-dimensional Nishimoto fractional derivatives, and Engheta \cite{Engheta1998}, based on one-dimensional Riemann-Liouville fractional derivatives.

\vspace*{-7pt}

\subsection{A fractional calculus by Tarasov}
We briefly describe the fractional calculus developed by Tarasov in Ref. \cite{Tarasov2008} and highlight its differences with respect to the unified nonlocal vector calculus that has been the focus of this article.
Tarasov's calculus is built upon one-dimensional Caputo derivatives (as opposed to Riemann-Liouville derivatives used in \eqref{eq:frac-div} and \eqref{eq:frac-grad}) and a fractional generalization of the Fundamental Theorem of Calculus. These are used to define fractional differential and integral vector operations which satisfy fractional Green's, Stokes' and Gauss's theorem.

Tarasov introduces the following definitions.
If $u(\xb)$ is a $(r-1)$ times continuously differentiable scalar field such that $\frac{\partial^{r-1}u}{\partial \xb_i}$ is absolutely continuous, then its {\it fractional gradient} is defined as
\vskip -10pt %
\begin{equation}\label{eq:gradTarasov}
\text{Grad}^s_W u(\xb)= \,
{\bf e}_1 \ ^C\!D^s_W[\xb_1] u(\xb)+ \,
{\bf e}_2 \ ^C\!D^s_W[\xb_2] u(\xb) + \,
{\bf e}_3 \ ^C\!D^s_W[\xb_3] u(\xb).
\end{equation}
If $\vb(\xb)$ is a $(r-1)$ times continuously differentiable vector field such that $\frac{\partial^{r-1}\vb_i}{\partial \xb_i}$ are absolutely continuous, then its {\it fractional divergence} is defined as
\vskip-12pt %
\begin{equation}\label{eq:divTarasov}
\text{Div}_W^s \vb(\xb)= \;
^C\!D^s_W[\xb_1] \vb_1(\xb)+ \;
^C\!D^s_W[\xb_2] \vb_2(\xb)+ \;
^C\!D^s_W[\xb_3] \vb_3(\xb),
\end{equation}
Here, $ ^C\!D^s_W[\xb_i]$ is the fractional Caputo derivative with respect to the $i$-th component of $\xb$ for the parallelepiped
$W:=[a,b]\!\times\![c,d]\!\times\![g,h]$, i.e.
$$
^C\!D^s_W [\xb_1]=\, _a^CD^s_b [\xb_1], \;\;\;
^C\!D^s_W [\xb_2]=\, _c^CD^s_d [\xb_2], \;\;\;
^C\!D^s_W [\xb_3]=\, _g^CD^s_h [\xb_3],
$$
and the one-dimensional Caputo derivative is defined as
\vskip -12pt
$$
_a^CD^s_b[x]f(x)=\frac{1}{\Gamma(r-s)}
\int^b_a  \frac{1}{(b-x)^{1+s-r}}
\frac{\partial^r f}{\partial {x}^r}(x) dx \quad (r-1<s<r) .
$$

An advantage of the vector calculus developed by Ref.~\cite{Tarasov2008} is that a fractional curl operator can also be defined analogously to \eqref{eq:gradTarasov} and \eqref{eq:divTarasov}. Then, a fractional Fundamental Theorem of Calculus relating the Caputo fractional derivative and the Riemann-Louiville fractional integral as inverses is proved and leveraged to establish fractional Green’s, Stokes’ and Gauss’s theorems. These results are then applied to propose and study a set of fractional Maxwell's equations. In contrast, Refs.~\cite{Adda1997} and \cite{Engheta1998} do not investigate such analogues of the classical theorems of vector calculus and Ref. \cite{Meerschaert2006} only considers hybrid$^{vi}$
\footnote{$^{vi}$
Hybrid refers to the fact that in Ref. \cite{Meerschaert2006} higher order fractional operators, such as the fractional Laplacian, are defined as a composition of a local (divergence) and a fractional (gradient) operator.
} versions of divergence and Stokes' theorems. Though relevant, we do not investigate such results in this work. However, we mention that the \emph{unweighted} nonlocal vector calculus provides such results \cite{Du2013}; the extension to the weighted nonlocal vector calculus used throughout our article is the subject of current research.

Since the formulas \eqref{eq:gradTarasov} and \eqref{eq:divTarasov} depend only of the values of the input in the neighborhood of the axes of $\mathbb{R}^d$ defined by the unit vectors $\mathbf{e}_i$, these operators are not equivalent to the operators $\text{grad}^s$ and $\text{div}^s$ defined by \eqref{eq:frac-grad} and \eqref{eq:frac-div}, respectively.
The latter depend on the values of their input over $\mathbb{R}^d$ if the measure $M$ is supported on $\mathbb{R}^d$.
A major difference between the two families of vector calculus operators is that the analogue of \eqref{eq:F-WN-equivalence} for the operators \eqref{eq:gradTarasov} and \eqref{eq:divTarasov} is
\vskip -13pt
\begin{multline}\label{tarasov_composition}
\text{Div}_W^s \text{Grad}_W^s  u(\xb)
\\ = \;
\left(^C\!D^{s}_W[\xb_1]\right)^2 u(\xb)+ \;
\left(^C\!D^{s}_W[\xb_2]\right)^2 u(\xb)+ \;
\left(^C\!D^{s}_W[\xb_3]\right)^2 u(\xb),
\end{multline}
which involves not the isotropic operator $(-\Delta)^s$ but rather an anisotropic fractional operator of order $2s$. While the fractional Laplacian $(-\Delta)^s$ is the (negative) generator of isotropic $\alpha$-stable L\'evy motion \cite{Meerschaert2006,Meerschaert2012}, it is unclear whether the operator appearing in \eqref{tarasov_composition} has any relation to a diffusion process.

\vspace*{-7pt}
\subsection{A fractional vector calculus based on invariance requirements}\label{sec:silhavy}

The recent article by {\v{S}}ilhav{\`y} \cite{silhavy2020fractional} provides a critique of coordinate-based approaches for fractional vector calculus, such as the framework discussed in the previous section. Ref.~\cite{silhavy2020fractional} points out that fractional gradient operators $G$ defined component-wise using one-dimensional fractional derivatives, such as  $G = \text{Grad}_W^s$ defined by \eqref{eq:gradTarasov}, do not satisfy the identity
\begin{equation}
G(u \circ \eta_\mathcal{R}) = \mathcal{R}(Gu)\circ \eta_{\mathcal{R}}
\end{equation}
for a rotation matrix $\mathcal{R}$ and $\eta_{\mathcal{R}}(\xb) = \mathcal{R}\xb$. In other words, such gradient operators do not transform under rotations in the same way as classical gradient operators \cite{spivak2018calculus}. Ref.~\cite{silhavy2020fractional} then proposes specific transformation rules for translation, rotation, and scaling (homogeneity) that are desired for a fractional vector calculus. Operators that are equivalent to \eqref{eq:frac-Cart}, up to constants, are introduced and shown to satisfy these requirements. Then, Ref.~\cite{silhavy2020fractional} shows that any operators which satisfy such transformation rules and a certain continuity assumption must be given by scalar multiples of the operators \eqref{eq:frac-Cart} for Schwartz test function input.

Notably, Ref.~\cite{silhavy2020fractional} considers extensions of the vector divergence and gradient \eqref{eq:frac-Cart} to complex orders using analytic continuation. Among several useful identities, the author points out that once these fractional vector operators are expressed in terms of Riesz potentials, the following identity follows from results of Horv{\'a}th \cite{horvath1959some}:
\vskip -10pt
\begin{equation}
\text{div}^\alpha \text{grad}^\beta = - C (-\Delta)^{\alpha+\beta},
\end{equation}
for a constant $C$, $\text{Re } \alpha \ge -n, \text{Re } \beta \ge 0$, and for the operators in \eqref{eq:frac-Cart}.
Note that with the equivalences established by Lemma \ref{lemma:polar-Cart-equivalence}, this furnishes an alternative proof of \eqref{thm:frac-div-grad-lapl} for Schwartz class functions. Unlike the work of Tarasov \cite{Tarasov2008}, the framework established by {\v{S}}ilhav{\`y} \cite{silhavy2020fractional} is equivalent to the unified calculus we consider.

\vspace*{-4pt}

\appendix
\renewcommand{\theequation}{\Alph{section}.\arabic{equation}} 

\section*{Appendix}

\section{Proof of Lemma \ref{lem:integral-equality-constant}} \label{sec:integral-equality-constant}

\setcounter{equation}{0} 

We begin analyzing
\begin{equation}\label{eq:fourier_integral}
\int_{|\thetab| = 1} \int_{|\thetab'| = 1}
\thetab \cdot \thetab' (i\thetab \cdot \xib)^{s} (i\thetab' \cdot \xib)^{s} d\thetab d\thetab'
\end{equation}
by following the same arguments of Example 6.24 in Ref. \cite{Meerschaert2012}. The key idea is to extract the dependence on $|\xib|$ from the integral and show, by symmetry, that the remaining factor does not depend on $\xib$.
First, note that
\vskip -15pt%
\begin{align*}
(i\xib \cdot \thetab)^s
&=
|\xib \cdot \thetab|^s
e^{i \text{sgn}(\xib \cdot \thetab) \pi s/2 } \\
&=
|\xib \cdot \thetab|^s
\left[
\cos(\pi s/2) + i \text{ sgn}(\xib \cdot \thetab)\sin(\pi s/2)
\right]
\end{align*}
\vskip -4pt \noindent%
and likewise
\vskip -12pt%
\begin{align*}
(i\xib \cdot \bm{\theta'})^s
&=
|\xib \cdot \bm{\theta'}|^s
\left[
\cos(\pi s/2) + i \text{ sgn}(\xib \cdot \bm{\theta'})\sin(\pi s/2)
\right]
\end{align*}
\vskip -4pt \noindent%
so that
\vskip -13pt
\begin{multline*}
(i\xib \cdot \thetab)^s
(i\xib \cdot \bm{\theta'})^s
 \\ =
|\xib \cdot \thetab|^s
|\xib \cdot \bm{\theta'}|^s
\Big(
\left[
\cos^2(\pi s/2) -
\text{ sgn}(\xib \cdot \thetab)
\text{ sgn}(\xib \cdot \bm{\theta'})
\sin^2(\pi s/2)
\right]
\\
+
i\cos(\pi s/2) \sin(\pi s/2)
\left[
\text{sgn}(\xib \cdot \bm{\theta'})
+
\text{sgn}(\xib \cdot \thetab)
\right]\Big).
\end{multline*}
Next, write $\xib = |\xib| \bm{\nu}$, where $|\bm{\nu}| = 1$. Then
\begin{multline*}
(i\xib \cdot \thetab)^s
(i\xib \cdot \bm{\theta'})^s
= \\
|\xib|^{2s}
|\bm{\nu} \cdot \thetab|^s
|\bm{\nu} \cdot \bm{\theta'}|^s
\Big(
\big[\cos^2(\pi s/2) -
\text{ sgn}(\bm{\nu} \cdot \thetab)
\text{ sgn}(\bm{\nu} \cdot \bm{\theta'})
\sin^2(\pi s/2)
\big] \\
+
i\cos(\pi s/2) \sin(\pi s/2)
\left[
\text{sgn}(\bm{\nu} \cdot \bm{\theta'})
+
\text{sgn}(\bm{\nu} \cdot \thetab)
\right]\Big).
\end{multline*}
Therefore, the imaginary part of the integral \eqref{eq:fourier_integral}
is proportional to the integral
\vskip -12pt
\begin{equation*}
\int_{|\thetab| = 1}
\int_{|\bm{\theta'}| = 1}
(\bm{\theta'} \cdot \thetab)
|\bm{\nu} \cdot \thetab|^s
|\bm{\nu} \cdot \bm{\theta'}|^s
\left[
\text{sgn}(\bm{\nu} \cdot \bm{\theta'})
+
\text{sgn}(\bm{\nu} \cdot \thetab)
\right]
d\bm{\theta'}
d\thetab,
\end{equation*}
which is zero by symmetry. More precisely, the change-of-variables
\begin{equation}
(\bm{\theta}, \bm{\theta}') \mapsto (-\bm{\theta}, -\bm{\theta}')
\end{equation}
shows the integral is the negative of itself. Thus, the integral \eqref{eq:fourier_integral} is real-valued and given by
\vskip -12pt%
\begin{multline}
|\xib|^{2s}
\int_{|\thetab| = 1}
\int_{|\bm{\theta'}| = 1}
|\bm{\nu} \cdot \thetab|^s
|\bm{\nu} \cdot \bm{\theta'}|^s
\big[\cos^2(\pi s/2) \\ -
\text{ sgn}(\bm{\nu} \cdot \thetab)
\text{ sgn}(\bm{\nu} \cdot \bm{\theta'})
\sin^2(\pi s/2)
\big]
(\bm{\theta'} \cdot \thetab)
d\bm{\theta'}
d\thetab.
\end{multline}
Next we argue that the integral above does not depend on $\xib$ (note that this is not obvious at first glance, because $\bm{\nu} = \xib/|\xib|$). Denote the integral in the above equation by $F(\bm{\nu})$. We seek to show that $F(\bm{\nu})$ is the constant $D_{n,s}$. Let $T$ be an orthonormal rotation; then
\vskip -12pt%
\begin{multline}
F(T\bm{\nu}) =
\int_{|\thetab| = 1}
\int_{|\bm{\theta'}| = 1}
|T\bm{\nu} \cdot \thetab|^s
|T\bm{\nu} \cdot \bm{\theta'}|^s
\big[\cos^2(\pi s/2)  \\ -
\text{ sgn}(T\bm{\nu} \cdot \thetab)
\text{ sgn}(T\bm{\nu} \cdot \bm{\theta'})
\sin^2(\pi s/2)
\big]
(\bm{\theta'} \cdot \thetab)
d\bm{\theta'}
d\thetab.
\end{multline}
The change of variable $\thetab \rightarrow T\thetab$, which has Jacobian determinant $1$ in absolute value, then gives
\vskip -13pt%
\begin{multline*}
F(T\bm{\nu}) =
\int_{|\thetab| = 1}
\int_{|\bm{\theta'}| = 1}
|T\bm{\nu} \cdot T\thetab|^s
|T\bm{\nu} \cdot T\bm{\theta'}|^s
\big[\cos^2(\pi s/2) \\ -
\text{ sgn}(T\bm{\nu} \cdot T\thetab)
\text{ sgn}(T\bm{\nu} \cdot T\bm{\theta'})
\sin^2(\pi s/2)
\big]
(T\bm{\theta'} \cdot T\thetab)
d\bm{\theta'}
d\thetab.
\end{multline*}
Since $T\mathbf{a} \cdot T\mathbf{b} = \mathbf{a} \cdot \mathbf{b}$,
\vskip -12pt%
\begin{align*}
F(T\bm{\nu}) &= \int_{|\thetab| = 1}
\begin{multlined}[t]
\int_{|\bm{\theta'}| = 1}
|\bm{\nu} \cdot \thetab|^s
|\bm{\nu} \cdot \bm{\theta'}|^s
\big[\cos^2(\pi s/2) \\ -
\text{ sgn}(\bm{\nu} \cdot \thetab)
\text{ sgn}(\bm{\nu} \cdot \bm{\theta'})
\sin^2(\pi s/2)
\big]
(\bm{\theta'} \cdot \thetab)
d\bm{\theta'}
d\thetab
\end{multlined}
\\
&= F(\nu).
\end{align*}
Therefore $F(\bm{\nu})$ does not depend on $\bm{\nu}$; for example $F(\bm{\nu}) = F(\mathbf{e}_1)$, by choosing an orthonormal $T$ such that $T(\bm{\nu}) = \mathbf{e}_1$. We can then write the above constant as
\vskip -13pt
\begin{multline}\label{eq:F_e_1}
F(\mathbf{e}_1) =
\cos^2(\pi s/2)
\int_{|\bm{\theta}| = 1}
\int_{|\bm{\theta'}| = 1}
|\bm{\theta}_1|^s
|\bm{\theta'}_1|^s
(\bm{\theta'} \cdot \bm{\theta})
d\bm{\theta'}
d\bm{\theta} \\
-
\sin^2(\pi s/2)
\int_{|\bm{\theta}| = 1}
\int_{|\bm{\theta'}| = 1}
|\bm{\theta}_1|^s
|\bm{\theta'}_1|^s
\text{ sgn}(\bm{\theta}_1)
\text{ sgn}(\bm{\theta'}_1)
(\bm{\theta'} \cdot \bm{\theta})
d\bm{\theta'}
d\bm{\theta}.
\end{multline}
\vskip -3pt \noindent
For the first term above, consider the transformation $\bm{\theta}' \mapsto -\bm{\theta'}$. Under this transformation, $\bm{\theta} \cdot \bm{\theta'} \mapsto -\bm{\theta} \cdot \bm{\theta'}$. This transformation also maps the unit sphere $\{|\bm{\theta}'|=1\}$ to itself and has Jacobian determinant $1$ in absolute value.
Therefore
\vskip -14pt
\begin{equation}
\begin{aligned}
&\cos^2(\pi s/2)
\int_{|\bm{\theta}| = 1}
\int_{|\bm{\theta'}| = 1}
|\bm{\theta}_1|^s
|\bm{\theta'}_1|^s
(\bm{\theta'} \cdot \bm{\theta})
d\bm{\theta'}
d\bm{\theta}\\
=&
-\cos^2(\pi s/2)
\int_{|\bm{\theta}| = 1}
\int_{|\bm{\theta'}| = 1}
|\bm{\theta}_1|^s
|\bm{\theta'}_1|^s
(\bm{\theta'} \cdot \bm{\theta})
d\bm{\theta'}
d\bm{\theta}
\end{aligned}
\end{equation}
\vskip -3pt \noindent
which proves the first term in \eqref{eq:F_e_1} is zero.
Therefore
\vskip- 13pt
\begin{equation}
D_{n,s} =
-
\sin^2(\pi s/2)
\int_{|\bm{\theta}| = 1}
\int_{|\bm{\theta'}| = 1}
|\bm{\theta}_1|^s
|\bm{\theta'}_1|^s
\text{ sgn}(\bm{\theta}_1)
\text{ sgn}(\bm{\theta'}_1)
(\bm{\theta'} \cdot \bm{\theta})
d\bm{\theta'}
d\bm{\theta}.
\end{equation}
\vskip -4pt \noindent
We write
\vskip -13pt
\begin{align*}
\int_{|\bm{\theta}| = 1}
\int_{|\bm{\theta'}| = 1}
&|\bm{\theta}_1|^s
|\bm{\theta'}_1|^s
\text{ sgn}(\bm{\theta}_1)
\text{ sgn}(\bm{\theta'}_1)
(\bm{\theta'} \cdot \bm{\theta})
d\bm{\theta'}
d\bm{\theta}
\\
&=
\int_{\substack{|\bm{\theta}| = 1 \\ \bm{\theta}_1 \ge 0}}
\int_{\substack{|\bm{\theta'}| = 1 \\ \bm{\theta'}_1 \ge 0}}
|\bm{\theta}_1|^s
|\bm{\theta'}_1|^s
\text{ sgn}(\bm{\theta}_1)
\text{ sgn}(\bm{\theta'}_1)
(\bm{\theta'} \cdot \bm{\theta})
d\bm{\theta'}
d\bm{\theta}
\\
&\quad+
\int_{\substack{|\bm{\theta}| = 1 \\ \bm{\theta}_1 \ge 0}}
\int_{\substack{|\bm{\theta'}| = 1 \\ \bm{\theta'}_1 \le 0}}
|\bm{\theta}_1|^s
|\bm{\theta'}_1|^s
\text{ sgn}(\bm{\theta}_1)
\text{ sgn}(\bm{\theta'}_1)
(\bm{\theta'} \cdot \bm{\theta})
d\bm{\theta'}
d\bm{\theta}
\\
&\quad+
\int_{\substack{|\bm{\theta}| = 1 \\ \bm{\theta}_1 \le 0}}
\int_{\substack{|\bm{\theta'}| = 1 \\ \bm{\theta'}_1 \le 0}}
|\bm{\theta}_1|^s
|\bm{\theta'}_1|^s
(\bm{\theta'} \cdot \bm{\theta})
d\bm{\theta'}
d\bm{\theta}
\\
&\quad+
\int_{\substack{|\bm{\theta}| = 1 \\ \bm{\theta}_1 \le 0}}
\int_{\substack{|\bm{\theta'}| = 1 \\ \bm{\theta'}_1 \ge 0}}
|\bm{\theta}_1|^s
|\bm{\theta'}_1|^s
\text{ sgn}(\bm{\theta}_1)
\text{ sgn}(\bm{\theta'}_1)
(\bm{\theta'} \cdot \bm{\theta})
d\bm{\theta'}
d\bm{\theta}
\end{align*}
\vskip -2pt \noindent
This can be simplified to
\begin{align*}
& \int_{|\bm{\theta}| = 1}
\int_{|\bm{\theta'}| = 1}
|\bm{\theta}_1|^s
|\bm{\theta'}_1|^s
\text{ sgn}(\bm{\theta}_1)
\text{ sgn}(\bm{\theta'}_1)
(\bm{\theta'} \cdot \bm{\theta})
d\bm{\theta'}
d\bm{\theta}
\\
&=
\int_{\substack{|\bm{\theta}| = 1 \\ \bm{\theta}_1 \ge 0}}
\int_{\substack{|\bm{\theta'}| = 1 \\ \bm{\theta'}_1 \ge 0}}
|\bm{\theta}_1|^s
|\bm{\theta'}_1|^s
(\bm{\theta'} \cdot \bm{\theta})
d\bm{\theta'}
d\bm{\theta}
-
\int_{\substack{|\bm{\theta}| = 1 \\ \bm{\theta}_1 \ge 0}}
\int_{\substack{|\bm{\theta'}| = 1 \\ \bm{\theta'}_1 \le 0}}
|\bm{\theta}_1|^s
|\bm{\theta'}_1|^s
(\bm{\theta'} \cdot \bm{\theta})
d\bm{\theta'}
d\bm{\theta}
\\
& +\,
\int_{\substack{|\bm{\theta}| = 1 \\ \bm{\theta}_1 \le 0}}
\int_{\substack{|\bm{\theta'}| = 1 \\ \bm{\theta'}_1 \le 0}}
|\bm{\theta}_1|^s
|\bm{\theta'}_1|^s
(\bm{\theta'} \cdot \bm{\theta})
d\bm{\theta'}
d\bm{\theta}
 -
\int_{\substack{|\bm{\theta}| = 1 \\ \bm{\theta}_1 \le 0}}
\int_{\substack{|\bm{\theta'}| = 1 \\ \bm{\theta'}_1 \ge 0}}
|\bm{\theta}_1|^s
|\bm{\theta'}_1|^s
(\bm{\theta'} \cdot \bm{\theta})
d\bm{\theta'}
d\bm{\theta}.
\end{align*}
In the third and fourth integrals, consider the change-of-variables
$(\bm{\theta}, \bm{\theta}') \mapsto -(\bm{\theta}, \bm{\theta}')$.
This maps the range of integration in the third term to that of the first term,
the range of integration in the fourth term to that of the second term, leaves
the integrands invariant, and has Jacobian determinant $1$ in absolute value.
Therefore, the above can be written
\vskip -10pt
\begin{equation}
2\int_{\substack{|\bm{\theta}| = 1 \\ \bm{\theta}_1 \ge 0}}
\int_{\substack{|\bm{\theta'}| = 1 \\ \bm{\theta'}_1 \ge 0}}
|\bm{\theta}_1|^s
|\bm{\theta'}_1|^s
(\bm{\theta'} \cdot \bm{\theta})
d\bm{\theta'}
d\bm{\theta}
-
2\int_{\substack{|\bm{\theta}| = 1 \\ \bm{\theta}_1 \ge 0}}
\int_{\substack{|\bm{\theta'}| = 1 \\ \bm{\theta'}_1 \le 0}}
|\bm{\theta}_1|^s
|\bm{\theta'}_1|^s
(\bm{\theta'} \cdot \bm{\theta})
d\bm{\theta'}
d\bm{\theta}.
\end{equation}
Finally, for the second term above, consider the transformation
$\bm{\theta}' \mapsto -\bm{\theta}'$. This maps the range of integration to
that of the first term, flips the sign of $\bm{\theta}' \cdot \bm{\theta}$, and
has Jacobian determinant $1$ in absolute value. The above difference of integrals can therefore be written
\begin{equation}\label{eq:one_integral}
4\int_{\substack{|\bm{\theta}| = 1 \\ \bm{\theta}_1 \ge 0}}
\int_{\substack{|\bm{\theta'}| = 1 \\ \bm{\theta'}_1 \ge 0}}
|\bm{\theta}_1|^s
|\bm{\theta'}_1|^s
(\bm{\theta'} \cdot \bm{\theta})
d\bm{\theta'}
d\bm{\theta}.
\end{equation}
Now write
\vskip -15pt
\begin{equation}
(\bm{\theta'} \cdot \bm{\theta})
=
\bm{\theta}_1 \bm{\theta'}_1
+
\sum_{i = 2}^n
\bm{\theta}_i \bm{\theta'}_i.
\end{equation}
Then we can write \eqref{eq:one_integral} as
\vskip -11pt
\begin{multline}\label{eq:sum_theta1_other}
4\int_{\substack{|\bm{\theta}| = 1 \\ \bm{\theta}_1 \ge 0}}
\int_{\substack{|\bm{\theta'}| = 1 \\ \bm{\theta'}_1 \ge 0}}
|\bm{\theta}_1|^s
|\bm{\theta'}_1|^s
\bm{\theta}_1 \bm{\theta'}_1
d\bm{\theta'}
d\bm{\theta} \\
+
4\int_{\substack{|\bm{\theta}| = 1 \\ \bm{\theta}_1 \ge 0}}
\int_{\substack{|\bm{\theta'}| = 1 \\ \bm{\theta'}_1 \ge 0}}
|\bm{\theta}_1|^s
|\bm{\theta'}_1|^s
\left(\sum_{i = 2}^n\bm{\theta}_i \bm{\theta'}_i\right)
d\bm{\theta'}
d\bm{\theta}.
\end{multline}
For the second integral here, consider the transformation
\vskip -10pt
\begin{equation}
(\bm{\theta}'_1, \bm{\theta}'_2, ..., \bm{\theta}'_n)
\mapsto
(\bm{\theta}'_1, -\bm{\theta}'_2, ..., -\bm{\theta}'_n).
\end{equation}
This maps the set $\{|\bm{\theta}'|=1, \bm{\theta}'_1 \ge 0\}$ to itself, since
it leaves both $|\bm{\theta}'| = (\bm{\theta}'_1)^2 + (\bm{\theta}'_2)^2 + ... (\bm{\theta}'_n)^2$ and
$\bm{\theta}'_1$ invariant. It also has Jacobian determinant $1$ in absolute value. But it flips the sign of the second integrand in \eqref{eq:sum_theta1_other}. Therefore, the second integral is zero. So, we are left with
\vskip -11pt
\begin{equation}
D_{n,s} =
-4\sin^2(\pi s/2)
\int_{\substack{|\bm{\theta}| = 1 \\ \bm{\theta}_1 \ge 0}}
\int_{\substack{|\bm{\theta'}| = 1 \\ \bm{\theta'}_1 \ge 0}}
|\bm{\theta}_1|^s
|\bm{\theta'}_1|^s
\bm{\theta}_1 \bm{\theta'}_1
d\bm{\theta'}
d\bm{\theta}.
\end{equation}
We see that the integral above is positive, since the integrand is positive over the region of integration. In fact, by Fubini's theorem, we have
\vskip -12pt
\begin{equation}
D_{n,s} =
-4\sin^2(\pi s/2)
\left(
\int_{\substack{|\bm{\theta}| = 1 \\ \bm{\theta}_1 \ge 0}}
|\bm{\theta}_1|^{s+1}
d\bm{\theta}
\right)^2.
\end{equation}

\vspace*{-4pt}

\section{The unweighted Green Identity \break for Fractional Operators}\label{sec:appendix_green}

\setcounter{equation}{0} 

In Section \ref{sec:frac-special-case}, we showed how \eqref{eq:unweighted-Green} for the fractional interaction kernel \eqref{eq:dipierro_alpha} yielded the unweighted fractional Green's identity \eqref{eq:unweighted-Green-fractional_MAIN}, i.e.,
\vskip -12pt %
\begin{multline*}
\frac{2}{C_{n,s}}
\int_\Omega (-\Delta)^s u(\xb) \, v(\xb) \,d\xb
=
\int_{\mathbb{R}^n} \int_{\mathbb{R}^n}
\frac{\left(u(\xb) - u(\yb)\right)\left(v(\xb) - v(\yb)\right)}{|\mathbf{x} - \mathbf{y}|^{n+2s}}
\,d\yb\,d\xb \\
-2 \int_{\mathbb{R}^n\setminus\Omega} \
v(\xb)
\int_{\mathbb{R}^n}
\frac{u(\mathbf{x}) - u(\mathbf{y})}
{|\mathbf{x} - \mathbf{y}|^{n+2s}}
\,d\yb\,d\xb.
\end{multline*}
The second term on the right-hand side can be split as
\vskip -12pt
\begin{align}
\label{second-term-rhs}
-2\int_{\mathbb{R}^n\setminus\Omega} \
v(\xb)
\int_{\mathbb{R}^n}
\frac{u(\mathbf{x}) - u(\mathbf{y})}
{|\mathbf{x} - \mathbf{y}|^{n+2s}}
\, &d\yb \,d\xb
\\
& = -2\int_{\mathbb{R}^n\setminus\Omega} \
v(\xb)
\int_{\mathbb{R}^n\setminus\Omega}
\frac{u(\mathbf{x}) - u(\mathbf{y})}
{|\mathbf{x} - \mathbf{y}|^{n+2s}}
\,d\yb \,d\xb
\\
&\ \ -2
\int_{\mathbb{R}^n\setminus\Omega} \
v(\xb)
\int_{\Omega}
\frac{u(\mathbf{x}) - u(\mathbf{y})}
{|\mathbf{x} - \mathbf{y}|^{n+2s}}
\,d\yb \,d\xb.
\end{align}
Meanwhile, the first term on the right-hand side can be split as:
\begin{align*}
\int\limits_{\mathbb{R}^n} \int\limits_{\mathbb{R}^n}
&\frac{\left(u(\xb) - u(\yb)\right)\left(v(\xb) - v(\yb)\right)}{|\mathbf{x} - \mathbf{y}|^{n+2s}}
\,d\yb\,d\xb
\\
&\qquad=
\int\limits_{\Omega} \int\limits_{\Omega}
\frac{\left(u(\xb) - u(\yb)\right)\left(v(\xb) - v(\yb)\right)}{|\mathbf{x} - \mathbf{y}|^{n+2s}}
\,d\yb\,d\xb \\
&\qquad+
\int\limits_{\mathbb{R}^n \setminus \Omega} \int\limits_{\Omega}
\frac{\left(u(\xb) - u(\yb)\right)\left(v(\xb) - v(\yb)\right)}{|\mathbf{x} - \mathbf{y}|^{n+2s}}
\,d\yb\,d\xb \\
&\qquad+
\int\limits_{\Omega} \int\limits_{\mathbb{R}^n \setminus \Omega}
\frac{\left(u(\xb) - u(\yb)\right)\left(v(\xb) - v(\yb)\right)}{|\mathbf{x} - \mathbf{y}|^{n+2s}}
\,d\yb\,d\xb \\
&\qquad+
\int\limits_{\mathbb{R}^n \setminus \Omega} \int\limits_{\mathbb{R}^n \setminus \Omega}
\frac{\left(u(\xb) - u(\yb)\right)\left(v(\xb) - v(\yb)\right)}{|\mathbf{x} - \mathbf{y}|^{n+2s}}
\,d\yb\,d\xb.
\end{align*}
\vskip -3pt \noindent %
The identity
\begin{align}
\begin{split}
(B \times B) \setminus (A \times A)
&=
\left((B \setminus A) \times (B \setminus A)\right) \cup
(A \times (B \setminus A)) \cup
((B \setminus A) \times A)
\end{split}
\end{align}
\vskip -5pt \noindent
implies, with $B = \mathbb{R}^n$, $A = \mathbb{R}^n \setminus \Omega$ and therefore $B \setminus A = \Omega$, that
\begin{equation}
\left(\mathbb{R}^n\right)^2 \setminus \left(\mathbb{R}^n \setminus \Omega\right)^2
=
(\Omega \times \Omega)
\cup
((\mathbb{R}^n \setminus \Omega) \times \Omega)
\cup
(\Omega \times (\mathbb{R}^n \setminus \Omega)).
\end{equation}
Therefore, the first three integrals can be combined into
\begin{equation}
\iint_{\left(\mathbb{R}^n\right)^2 \setminus \left(\mathbb{R}^n \setminus \Omega\right)^2}
\frac{\left(u(\xb) - u(\yb)\right)\left(v(\xb) - v(\yb)\right)}{|\mathbf{x} - \mathbf{y}|^{n+2s}}
\,d\yb\,d\xb.
\end{equation}
The fourth term can be written as
\begin{align*}
&\int_{\mathbb{R}^n \setminus \Omega} \int_{\mathbb{R}^n \setminus \Omega}
\frac{\left(u(\xb) - u(\yb)\right)\left(v(\xb) - v(\yb)\right)}{|\mathbf{x} - \mathbf{y}|^{n+2s}}
\,d\yb\,d\xb \\
&=
\begin{multlined}[t]
\int_{\mathbb{R}^n \setminus \Omega} \int_{\mathbb{R}^n \setminus \Omega}
v(\xb)
\frac{\left(u(\xb) - u(\yb)\right)}
{|\mathbf{x} - \mathbf{y}|^{n+2s}}
\,d\yb\,d\xb
\\
-
\int_{\mathbb{R}^n \setminus \Omega} \int_{\mathbb{R}^n \setminus \Omega}
v(\yb)
\frac{\left(u(\xb) - u(\yb)\right)}{|\mathbf{x} - \mathbf{y}|^{n+2s}}
\,d\yb\,d\xb
\end{multlined}
\\
&=
\begin{multlined}[t]
\int_{\mathbb{R}^n \setminus \Omega} v(\xb)
\int_{\mathbb{R}^n \setminus \Omega}
\frac{\left(u(\xb) - u(\yb)\right)}
{|\mathbf{x} - \mathbf{y}|^{n+2s}}
\,d\yb\,d\xb
\\
+
\int_{\mathbb{R}^n \setminus \Omega}
v(\xb)
\int_{\mathbb{R}^n \setminus \Omega}
\frac{\left(u(\xb) - u(\yb)\right)}{|\mathbf{x} - \mathbf{y}|^{n+2s}}
\,d\yb\,d\xb
\end{multlined}
\\
&=
2\int_{\mathbb{R}^n \setminus \Omega}
v(\xb)
\int_{\mathbb{R}^n \setminus \Omega}
\frac{\left(u(\xb) - u(\yb)\right)}{|\mathbf{x} - \mathbf{y}|^{n+2s}}
\,d\yb\,d\xb.
\end{align*}
This cancels with the first term in \eqref{second-term-rhs}, so \eqref{eq:unweighted-Green-fractional_MAIN} reduces to
\begin{align}\label{eq:dipierro-equivalent}
\frac{2}{C_{n,s}}
\int_\Omega (-\Delta)^s u(\xb) \, &v(\xb) \,d\xb
\\
&= \int \int_{\left(\mathbb{R}^n\right)^2 \setminus \left(\mathbb{R}^n \setminus \Omega\right)^2} \hspace{-1mm}
\frac{\left(u(\xb) - u(\yb)\right)\left(v(\xb) - v(\yb)\right)}{|\mathbf{x} - \mathbf{y}|^{n+2s}}
\,d\yb\,d\xb \\
&\qquad \qquad \qquad
-2 \int_{\mathbb{R}^n\setminus\Omega} \
v(\xb)
\int_{\Omega}
\frac{u(\mathbf{x}) - u(\mathbf{y})}
{|\mathbf{x} - \mathbf{y}|^{n+2s}}
\,d\yb\,d\xb.
\end{align}

This is precisely the fractional Green's identity \eqref{eq:dipierro_MAIN} reported by Ref.~\cite{Dipierro2017}. Thus, the latter can be viewed as a special case of the unweighted nonlocal Green's identity of Ref.~\cite{Du2013} with the fractional interaction kernel \eqref{eq:dipierro_alpha}.

Note that Ref.~\cite{Du2013} also defines a nonlocal interaction operator for a two-point vector field $\wb(\xb,\yb)$ and for $\xb\in\mbRn\setminus\Omega$,
\begin{equation}
\mathcal{N} \left[ \wb \right](\xb) =
- \int_{\mathbb{R}^n} (\wb(\xb,\yb) + \wb(\yb,\xb)) \cdot \alphab(\xb,\yb) d\yb
\end{equation}
which they use to introduce a notion of nonlocal flux,
\begin{align}
\mathcal{N} \big[ \mcG u \big] (\xb) &=
- \int_{\mathbb{R}^n} \big(\mcG u (\xb,\yb) + \mcG u (\yb,\xb) \big) \cdot \alphab(\xb,\yb) d\yb.
\end{align}
With the same fractional interaction kernel $\alphab$ given by \eqref{eq:dipierro_alpha}, using \eqref{eq:unw-grad} we have
\vskip -13pt %
\begin{align*}
\mathcal{N} \big[ \mcG u \big] (\xb)
= 2 \int_{\mathbb{R}^n}
\frac{u(\xb)-u(\yb)}{|\mathbf{x} - \mathbf{y}|^{n+2s}} d\yb.
\end{align*}
\vskip -3pt \noindent
The identity \eqref{eq:unweighted-Green-fractional_MAIN} can therefore be rewritten using the notation $\mathcal{N} \big[ \mcG u \big] (\xb)$, just as \eqref{eq:dipierro_MAIN} uses the notation $\mathcal{N}_s u(\xb)$. The two operators are not equivalent, as
\vskip -14pt
\begin{equation}
\mathcal{N} \big[ \mcG u \big] (\xb) =
2\mathcal{N}_s u(\xb)
+
2 \int_{\mathbb{R}^n \setminus \Omega}
\frac{u(\xb)-u(\yb)}{|\mathbf{x} - \mathbf{y}|^{n+2s}} d\yb.
\end{equation}


\vspace*{-7pt}

\section{Proof of Theorem \ref{thm:L2convergence}} \label{sec:L2convergence}

\setcounter{equation}{0} 

Theorem \ref{thm:L2convergence} is a direct consequence of the following two lemmas.

\begin{manuallemma}{C.1}
\label{lem:grad}
{\sl Let $u \in H^s(\mbRn)$ and $\Upsilon_{n,t}$ be defined as in \eqref{eq:Upsilon}. Then
\vskip -12pt %
\begin{multline}\label{eq:grad-lemma}
      \int_{\mbRn} \left|\int_{\mbRn \setminus B_\delta(\xb)}\right.
     \left.\left[ u(\xb) - u(\yb) \right]
      \frac{\xb-\yb}{|\xb-\yb|}
      \frac{1}{|\xb-\yb|^{n+s}}d\yb \right|^2 d\xb \\
\le \| u \|_{L_2(\mbRn)}^2
      \left(\frac{2\Upsilon_{n,n+2s}}{\delta^{n+2s}}
      + \|u\|_{L^1(\mbRn)}
      \sqrt{\frac{\Upsilon_{n,3n+4s}}{\delta^{3n+4s}}}\right).
\end{multline}
}
\end{manuallemma}

\proof
We find a bound for the left-hand side in \eqref{eq:grad-lemma}.
First,
\vskip -13pt %
\begin{align*}
     \Bigg| \int_{\mbRn \setminus B_\delta(\xb)}
     \left[ u(\xb) - u(\yb) \right]&
     \frac{\xb-\yb}{|\xb-\yb|}
     \frac{1}{|\xb-\yb|^{n+s}}d\yb \Bigg|^2 \\
&\le \int_{\mbRn \setminus B_\delta(\xb)}
     \left| \left[ u(\xb) - u(\yb) \right]
     \frac{\xb-\yb}{|\xb-\yb|}
     \frac{1}{|\xb-\yb|^{n+s}}\right|^2 d\yb \\
&\le \int_{\mbRn \setminus B_\delta(\xb)}
     \frac{(u(\xb) - u(\yb))^2}{|\xb-\yb|^{2(n+s)}}d\yb\\
&=   \int_{\mbRn \setminus B_\delta(\xb)}
     \frac{(u(\xb))^2 -2u(\xb)u(\yb)+(u(\yb))^2}
     {|\xb-\yb|^{2(n+s)}}d\yb.
\end{align*}
\vskip -3pt \noindent %
Then
\vskip -14pt %
$$
\int_{\mbRn}  \Bigg|\int_{\mbRn \setminus B_\delta(\xb)}
\left[ u(\xb) - u(\yb) \right] \frac{\xb-\yb}{|\xb-\yb|}
\frac{1}{|\xb-\yb|^{n+s}}d\yb
\Bigg|^2 d\xb
$$
\begin{align*}
&\begin{multlined}[t]
\le \int_{\mbRn} \int_{\mbRn \setminus B_\delta(\xb)}
\frac{(u(\xb))^2}{|\xb-\yb|^{2(n+s)}} d\yb d\xb
- 2 \int_{\mbRn}\int_{\mbRn \setminus B_\delta(\xb)}
\frac{u(\xb)u(\yb)}{|\xb-\yb|^{2(n+s)}} d\yb d\xb \\
+ \int_{\mbRn}\int_{\mbRn \setminus B_\delta(\xb)}
\frac{(u(\yb))^2}{|\xb-\yb|^{2(n+s)}} d\yb d\xb
\end{multlined}\\
&= I + II + III.
\end{align*}
We have from Lemma \ref{upsilon_lemma} with $t = n+2s$,
\begin{equation*}
I = \int_{\mbRn}(u(\xb))^2
\int_{\mbRn \setminus B_\delta(\xb)}
\frac{1}{|\xb-\yb|^{n+(n+2s)}}d\yb d\xb
=\frac{\Upsilon_{n,n+2s}}{\delta^{n+2s}}\| u \|_{L_2(\mbRn)}^2.
\end{equation*}
Next, by using Fubini's theorem, we see that $III = I$; in fact, the change of variables $(\xb,\yb) \mapsto (\yb,\xb)$, and the symmetry of $|\xb-\yb|$ yields
\vskip -11pt%
\begin{align*}
III &= \int_{\mbRn}
\int_{\mbRn \setminus B_\delta(\xb)}
\frac{(u(\yb))^2}{|\xb-\yb|^{2(n+s)}} d\yb d\xb \\
&=
\iint_{\{ (\xb,\yb) \ : \ |\xb-\yb| \ge \delta\}}
\frac{(u(\yb))^2}{|\xb-\yb|^{2(n+s)}} d\yb d\xb \\
&=
\iint_{\{ (\xb,\yb) \ : \ |\xb-\yb| \ge \delta\}}
\frac{(u(\xb))^2}{|\xb-\yb|^{2(n+s)}} d\yb d\xb \\
&= \int_{\mbRn}\int_{\mbRn \setminus B_\delta(\xb)}
\frac{(u(\xb))^2}{|\yb-\xb|^{2(n+s)}} d\yb d\xb = I.
\end{align*}
To bound $II$, we first estimate, using Lemma \ref{upsilon_lemma} with $t = 3n+4s$,
\vskip -10pt%
\begin{align*}
|u(\xb)|
\Bigg|\int_{\mbRn \setminus B_\delta(\xb)}
&\frac{u(\yb)}{|\xb-\yb|^{2(n+s)}} d\yb\Bigg|
\le
|u(\xb)|\left|\int_{\mbRn \setminus B_\delta(\xb)}
\frac{u(\yb)}{|\xb-\yb|^{2(n+s)}}d\yb\right| \\
&\qquad \le |u(\xb)| \| u \|_{L^2(\mbRn)}
\sqrt{\int_{\mbRn \setminus B_\delta(\xb)}
\frac{1}{|\xb-\yb|^{4(n+s)}}d\yb} \\
&\qquad=
|u(\xb)| \| u \|_{L^2(\mbRn)}
\sqrt{\int_{\mbRn \setminus B_\delta(\xb)}
\frac{1}{|\xb-\yb|^{n+(3n+4s))}}d\yb} \\
&\qquad\le|u(\xb)| \| u \|_{L^2(\mbRn)}
\sqrt{\frac{\Upsilon_{n,3n+4s}}{\delta^{3n+4s}}},
\end{align*}
so that
\begin{align*}
\int_{\mbRn}
|u(\xb)|
\left|\int_{\mbRn \setminus B_\delta(\xb)}
\frac{u(\yb)}{|\xb-\yb|^{2(n+s)}}d\yb\right|d\xb
&\le\int |u(\xb)| \| u \|_{L^2(\mbRn)}
\sqrt{\frac{\Upsilon_{n,3n+4s}}{\delta^{3n+4s}}}d\xb \\
&=\|u\|_{L^1(\mbRn)} \| u \|_{L^2(\mbRn)}
\sqrt{\frac{\Upsilon_{n,3n+4s}}{\delta^{3n+4s}}}.
\end{align*}
Combining this bound for $II$ with the bound for $I = III$ yields \eqref{eq:grad-lemma}.
\endproof

\begin{manuallemma}{C.2}
\label{lem:div}
{\sl
Let $\vb \in \mathbf{H}^s(\mbRn)$ and $\Upsilon_{n,t}$ be defined as in \eqref{eq:Upsilon}. Then
\vskip -10pt %
\begin{multline}\label{eq:div-lemma}
\int_{\mbRn} \left|\int_{\mbRn \setminus B_\delta(\xb)}
\left[ \vb(\xb) - \vb(\yb) \right]\right.
\left.\cdot \frac{\xb-\yb}{|\xb-\yb|}
\frac{1}{|\xb-\yb|^{n+s}}d\yb\right|^2d\xb \\
\le \| \ |\vb| \ \|_{L_2(\mbRn)}^2 \left(
\frac{2\Upsilon_{n,n+2s}}{\delta^{n+2s}}
+\| \ |\vb| \ \|_{L^1(\mbRn)}
\sqrt{\frac{\Upsilon_{n,3n+4s}}{\delta^{3n+4s}}}\right).
\end{multline}
}
\end{manuallemma}

\proof
We proceed as before. First, we write
\begin{align*}
     \Bigg| \int_{\mbRn \setminus B_\delta(\xb)}
     \left[ \vb(\xb) - \vb(\yb) \right] &
     \cdot \frac{\xb-\yb}{|\xb-\yb|}
     \frac{1}{|\xb-\yb|^{n+s}}d\yb \Bigg|^2 \\
&\le \int_{\mbRn \setminus B_\delta(\xb)}
     \left| \left[ \vb(\xb) - \vb(\yb) \right]
     \frac{\xb-\yb}{|\xb-\yb|}
     \frac{1}{|\xb-\yb|^{n+s}}\right|^2 d\yb \\
&\le \int_{\mbRn \setminus B_\delta(\xb)}
     \frac{| \vb(\xb) - \vb(\yb)|^2}{|\xb-\yb|^{2(n+s)}}d\yb\\
&=   \int_{\mbRn \setminus B_\delta(\xb)}
     \frac{|\vb(\xb)|^2 -2\vb(\xb) \cdot \vb(\yb)+|\vb(\yb)|^2}
     {|\xb-\yb|^{2(n+s)}}d\yb.
\end{align*}
\vskip -2pt \noindent %
Then
\vskip -12pt %
\begin{align*}
\int_{\mbRn} & \Bigg|\int_{\mbRn \setminus B_\delta(\xb)}
\left[ \vb(\xb) - \vb(\yb) \right] \cdot \frac{\xb-\yb}{|\xb-\yb|}
\frac{1}{|\xb-\yb|^{n+s}}d\yb
\Bigg|^2 d\xb \\
&\begin{multlined}[t]
\le \int_{\mbRn} \int_{\mbRn \setminus B_\delta(\xb)}
\frac{|\vb(\xb)|^2}{|\xb-\yb|^{2(n+s)}} d\yb d\xb
- 2 \int_{\mbRn}\int_{\mbRn \setminus B_\delta(\xb)}
\frac{\vb(\xb) \cdot \vb(\yb)}{|\xb-\yb|^{2(n+s)}} d\yb d\xb \\
+ \int_{\mbRn}\int_{\mbRn \setminus B_\delta(\xb)}
\frac{|\vb(\yb)|^2}{|\xb-\yb|^{2(n+s)}} d\yb d\xb.
\end{multlined}\\
\end{align*}
\vskip -4pt \noindent
These terms can be compared to $I, II,$ and $III$ in the proof of Lemma \ref{lem:grad}, and the proof can completed in the same way.
\endproof

\vspace*{-6pt}

\section{Proof of Theorem \ref{thm:point-convergence}} \label{sec:point-convergence}

\setcounter{equation}{0} 

Theorem \ref{thm:point-convergence} is a direct consequence of the following two lemmas.

\begin{manuallemma}{D.1}
\label{lem:grad-pointwise}
{\sl
Let $u \in H^s(\mbRn)$ and $\Upsilon_{n,t}$ be defined as in \eqref{eq:Upsilon}. Then
\vskip -12pt %
\begin{multline}
\left|\;\int\limits_{\mbRn \setminus B_\delta(\xb)}
\left[ u(\xb) - u(\yb) \right]\frac{\xb-\yb}{|\xb-\yb|}
\frac{1}{|\xb-\yb|^{n+s}}d\yb\,\right|
\\
\le |u(\xb)|\frac{\Upsilon_{n,s}}{\delta^s} +
\|u\|_{L^2(\mbRn)}\sqrt{
\frac{\Upsilon_{n,n+2s}}{\delta^{n+2s}}}.
\end{multline}
}
\end{manuallemma}

\proof
A simple calculation gives
\vskip -13pt
\begin{align*}
\Bigg| \int_{\mbRn \setminus B_\delta(\xb)}
&\left[ u(\xb) - u(\yb) \right]
\frac{\xb-\yb}{|\xb-\yb|}
\frac{1}{|\xb-\yb|^{n+s}} d\yb\, \Bigg| \\
&\le
\int_{\mbRn \setminus B_\delta(\xb)}
\left[ |u(\xb)| + |u(\yb)| \right]
\frac{1}{|\xb-\yb|^{n+s}}
d\yb \\
&\le
|u(\xb)|
\int_{\mbRn \setminus B_\delta(\xb)}
\frac{1}{|\xb-\yb|^{n+s}}
d\yb
+
\int_{\mbRn \setminus B_\delta(\xb)}
|u(\yb)|
\frac{1}{|\xb-\yb|^{n+s}}
d\yb \\
&\le
|u(\xb)|
\frac{\Upsilon_{n,s}}{\delta^s} +
\|u\|_{L^2(\mbRn)}
\sqrt{
\int_{\mbRn \setminus B_\delta(\xb)}
\frac{1}{|\xb-\yb|^{2(n+s)}}
d\yb} \\
&=
|u(\xb)|
\frac{\Upsilon_{n,s}}{\delta^s} +
\|u\|_{L^2(\mbRn)}
\sqrt{
\int_{\mbRn \setminus B_\delta(\xb)}
\frac{1}{|\xb-\yb|^{n+(n+2s)}}
d\yb} \\	
&\le
|u(\xb)|
\frac{\Upsilon_{n,s}}{\delta^s} +
\|u\|_{L^2(\mbRn)}
\sqrt{
\frac{\Upsilon_{n,n+2s}}{\delta^{n+2s}}}.
\end{align*}
Here, we used Lemma \ref{upsilon_lemma} with $t = n+2s$.
\endproof

\begin{manuallemma}{D.2}
\label{lem:div-pointwise}
{\sl
Let $\vb \in {\bf H}^s(\mbRn)$ and $\Upsilon_{n,t}$ be defined as in \eqref{eq:Upsilon}. Then
\vskip - 12pt
\begin{multline}
\left|
\int_{\mbRn \setminus B_\delta(\xb)}
\left[ \vb(\xb) - \vb(\yb) \right]
\cdot
\frac{\xb-\yb}{|\xb-\yb|}
\frac{1}{|\xb-\yb|^{n+s}}
d\yb
\right| \\
\le |\vb(\xb)|
\frac{\Upsilon_{n,s}}{\delta^s} +
\| \ |\vb| \ \|_{L^2(\mbRn)}
\sqrt{
\frac{\Upsilon_{n,n+2s}}{\delta^{n+2s}}}.
\end{multline}
}
\end{manuallemma}

\proof
Following similar arguments as in the proof of Lemma \ref{lem:grad-pointwise}, we have
\vskip -12pt %
\begin{align*}
     \Bigg|\int_{\mbRn \setminus B_\delta(\xb)}
&    \left[ \vb(\xb) - \vb(\yb) \right]
     \cdot\frac{\xb-\yb}{|\xb-\yb|}
     \frac{1}{|\xb-\yb|^{n+s}} d\yb \,\Bigg|\\
&=   \left|\int_{\mbRn \setminus B_\delta(\xb)}
     \left[ |\vb(\xb)| + |\vb(\yb)| \right]
     \frac{1}{|\xb-\yb|^{n+s}}d\yb \right| \\
&\le |\vb(\xb)|\int_{\mbRn \setminus B_\delta(\xb)}
     \frac{1}{|\xb-\yb|^{n+s}}d\yb+
     \int_{\mbRn \setminus B_\delta(\xb)}
     |\vb(\yb)|\frac{1}{|\xb-\yb|^{n+s}}d\yb .
\end{align*}
vskip -2pt \noindent
The proof can now be completed following the same steps as in the proof of Lemma \ref{lem:grad-pointwise}.
\endproof

\vspace*{-14pt}


\section{Proof of Theorem \ref{thm:tempered-equivalence-kernel}}
\label{sec:proof_tempered_equivalence_kernel}

\setcounter{equation}{0} 

With the choices in Theorem \ref{thm:tempered-equivalence-kernel}, since the kernel is translation invariant, we write the equivalence kernel according to \eqref{eq:kerndef-transl}. We have
\vskip - 12pt%
\begin{equation*}
2\gameq(\xb,\yb)=
\int_\mbRn \frac{\yb-\zb}{|\yb-\zb|^{n+s+1}}\cdot
   \frac{\zb-\xb}{|\zb-\xb|^{n+s+1}}
e^{-\lambda|\yb-\zb|} e^{-\lambda|\zb-\xb|} d\zb.
\end{equation*}
We explicitly denote the dependence of $\gameq(\xb,\yb)$ on $\lambda$ by writing $\gameq(\xb,\yb; \lambda)$. We evaluate this integral indirectly. Let $\zb' = \zb-\xb$. Then $\zb = \zb' + \xb$, $d\zb = d\zb'$ and $\yb - \zb = \yb - \zb' - \xb = \yb - \xb - \zb'$. Thus,
\vskip -12pt %
\begin{align*}
2\gameq(\xb,\yb; \lambda) &= \int_{\mbRn}\frac{\yb-\xb-\zb'}
{|\yb-\xb-\zb'|^{n+s+1}}\cdot
\frac{\zb'}{|\zb'|^{n+s+1}}
e^{-\lambda|\yb-\xb-\zb'|} e^{-\lambda|\zb'|}
d\zb'\\
&=\int_{\mbRn}\frac{\yb-\xb-\zb}{|\yb-\xb-\zb|^{n+s+1}}
\cdot\frac{\zb}{|\zb|^{n+s+1}}
e^{-\lambda|\yb-\xb-\zb|} e^{-\lambda|\zb|}
d\zb.
\end{align*}
From this, it follows that $\gameq(\xb,\yb; \lambda)$ depends only on $\xb-\yb$. Next, we show that the kernel is rotationally invariant. Consider a rotation $\mathcal{R}$; we have
\vskip -13pt %
\begin{multline*}
2\gameq\left(\mathcal{R}(\xb-\yb); \lambda \right) \\
= \int_{\mbRn}
\frac{\mathcal{R}(\yb-\xb)-\zb}{|\mathcal{R}(\yb-\xb)-\zb|^{n+s+1}}
\cdot \frac{\zb}{|\zb|^{n+s+1}}
e^{-\lambda|\mathcal{R}({\yb-\xb})-\zb|} e^{-\lambda|\zb|}
d\zb.
\end{multline*}
Let $\zb = \mathcal{R} \zb'$. Then $d\zb = d\zb'$, and
\vskip -12pt %
\begin{align*}
2\gameq&\left(\mathcal{R}(\xb-\yb); \lambda \right) \\ &=
\int_{\mbRn}\frac{\mathcal{R}(\yb-\xb)-\mathcal{R}\zb'}
{|\mathcal{R}(\yb-\xb)-\mathcal{R}\zb'|^{n+s+1}}
\cdot\frac{\mathcal{R}\zb'}
{|\mathcal{R}\zb'|^{n+s+1}}
e^{-\lambda|\mathcal{R}({\yb-\xb})-\mathcal{R}\zb'|}
e^{-\lambda|\mathcal{R}\zb'|}
d\zb'
\\
&=
\int_{\mbRn}\frac{\mathcal{R}(\yb-\xb)-\mathcal{R}\zb}
{|\mathcal{R}(\yb-\xb)-\mathcal{R}\zb|^{n+s+1}}
\cdot\frac{\mathcal{R}\zb}{|\mathcal{R}\zb|^{n+s+1}}
e^{-\lambda|\mathcal{R}({\yb-\xb})-\mathcal{R}\zb|}
e^{-\lambda|\mathcal{R}\zb|}
d\zb
\end{align*} \begin{align*} 
&=
\int_{\mbRn}\frac{\mathcal{R}\left((\yb-\xb)-\zb\right)}
{|\mathcal{R}\left((\yb-\xb)-\zb\right)|^{n+s+1}}\cdot
\frac{\mathcal{R}\zb}{|\mathcal{R}\zb|^{n+s+1}}
e^{-\lambda|\mathcal{R}(({\yb-\xb})-\zb)|}
e^{-\lambda|\mathcal{R}\zb|}
d\zb
\\
&=
\begin{multlined}[t]
\int_{\mbRn}
\frac{1}{|\mathcal{R}\left((\yb-\xb)-\zb\right)|^{n+s+1}}
\frac{1}{|\mathcal{R}\zb|^{n+s+1}}
\left[{\mathcal{R}\left((\yb-\xb)-\zb\right)}\cdot
{\mathcal{R}\zb}\right] \\ \times
e^{-\lambda|\mathcal{R}(({\yb-\xb})-\zb)|}
e^{-\lambda|\mathcal{R}\zb|}
d\zb
\end{multlined}
\\
&=
\begin{multlined}[t]
\int_{\mbRn}
\frac{1}{|\left((\yb-\xb)-\zb\right)|^{n+s+1}}\frac{1}{|\zb|^{n+s+1}}
\left[{\left((\yb-\xb)-\zb\right)}\cdot{\zb}\right]
\\
\times
e^{-\lambda|({\yb-\xb})-\zb|}
e^{-\lambda|\zb|}
d\zb
\end{multlined}
\\
&=
\int_{\mbRn}
\frac{(\yb-\xb)-\zb}{|(\yb-\xb)-\zb|^{n+s+1}}
\cdot\frac{\zb}{|\zb|^{n+s+1}}
e^{-\lambda|({\yb-\xb})-\zb|}
e^{-\lambda|\zb|}
d\zb .
\end{align*}
Therefore, $\gameq$ depends only on $|\xb-\yb|$. Now we let $c > 0$ and consider
\vskip -13pt %
\begin{align*}
2\gameq&(c|\xb-\yb|; \lambda) =
\int_{\mbRn} \frac{c(\yb-\xb)-\zb}
{|c(\yb-\xb)-\zb|^{n+s+1}}
\cdot \frac{\zb}{|\zb|^{n+s+1}}
e^{-\lambda|c({\yb-\xb})-\zb|}
e^{-\lambda|\zb|}
d\zb.
\end{align*}
Let $\zb = c \zb'$. Then $d\zb = c^{d} d\zb'$, and
\vskip - 12pt
\begin{align*}
&2\gameq(c|\xb-\yb|; \lambda) \\
&=
\int_{\mbRn}\frac{c(\yb-\xb)-c \zb'}
{|c(\yb-\xb)-c \zb'|^{n+s+1}}
\cdot\frac{c \zb'}{|c \zb'|^{n+s+1}}
e^{-\lambda|c({\yb-\xb})-c\zb'|}
e^{-\lambda|c\zb'|}
c^n d\zb' \\
&=
\int_{\mbRn}
\frac{c(\yb-\xb)-c \zb}
{|c(\yb-\xb)-c \zb|^{n+s+1}}
\cdot\frac{c \zb}{|c \zb|^{n+s+1}}
e^{-\lambda|c({\yb-\xb})-c\zb|}
e^{-\lambda|c\zb|}
c^n d\zb \\
&=
\frac{c}{c^{n+s+1}}\frac{c}{c^{n+s+1}}c^n
\int_{\mbRn}\frac{(\yb-\xb)-\zb}{|(\yb-\xb)-\zb|^{n+s+1}}
\cdot\frac{\zb}{|\zb|^{n+s+1}}
e^{-c\lambda|({\yb-\xb})-\zb|}
e^{-c\lambda|\zb|}
d\zb \\
&=
\frac{1}{c^{n+2s}}
\int_{\mbRn}\frac{(\yb-\xb)-\zb}{|(\yb-\xb)-\zb|^{n+s+1}}\cdot
\frac{\zb}{|\zb|^{n+s+1}}
e^{-c\lambda|({\yb-\xb})-\zb|}
e^{-c\lambda|\zb|}
d\zb \\
&=
\frac{1}{c^{n+2s}}2\gameq(|\xb-\yb|; c\lambda).
\end{align*}
Therefore,
\vskip -13pt
\begin{align*}
2\gameq(\xb,\yb; \lambda) &=
2\gameq\left( |\xb - \yb| \frac{\xb - \yb}{|\xb - \yb|} ; \lambda \right)
= 2\gameq\left( |\xb - \yb| \bf{e} ; \lambda \right) \\
&= \frac{1}{|\xb-\yb|^{n+2s}}
2\gameq\left({\bf e}; {\lambda}{|\xb - \yb|}\right),
\end{align*}
\vskip -2pt \noindent %
where $\bf e$ is any unit vector and
\vskip -10pt %
\begin{equation*}
2\gameq({\bf e}; {\lambda}{|\xb - \yb|})
=
\int_{\mbRn}\frac{{\bf e}-\zb}{|{\bf e}-\zb|^{n+s+1}}
\cdot \frac{\zb}{|\zb|^{n+s+1}}
e^{-\lambda|\xb - \yb||{\bf e}-\zb|}
e^{-\lambda|\xb - \yb||\zb|}
d\zb.
\end{equation*}

\vspace*{-3pt}

\section*{Acknowledgments}

This work was supported by the U.S. Department of Energy, Office of Science, Office of Advanced Scientific Computing Research under the Collaboratory on Mathematics and Physics-Informed Learning Machines for Multiscale and Multiphysics Problems (PhILMs) project, as well as MURI/ARO grant W911NF-15-1-0562.

M.D. and M.G. are supported by Sandia National Laboratories (SNL), SNL is a multimission laboratory managed and operated by National Technology and Engineering Solutions of Sandia, LLC., a wholly owned subsidiary of Honeywell International, Inc., for the U.S. Department of Energy's National Nuclear Security Administration contract number \\ {DE-NA0003525}. This paper describes objective technical results and analysis. Any subjective views or opinions that might be expressed in the paper do not necessarily represent the views of the U.S. Department of Energy or the United States Government. Report Number: SAND2020-4869.

\smallskip

The authors are thankful to Prof. Mark M. Meerschaert$\dagger$ (Michigan State University) for useful discussions and insights and for having initiated and inspired this work at its preliminary stage.

\vspace*{-7pt}

\renewcommand{\bibliofont}{\normalsize}


\medskip   \medskip

 \it

 \noindent
$^1$ Computational Science and Analysis \\
Sandia National Laboratories \\
Livermore, CA, 94550, USA  \\[3pt]
  e-mail: mdelia@sandia.edu \ (Corresponding author)\\
\\[3pt]
$^2$ Computational Science and Analysis \\
Sandia National Laboratories \\
Livermore, CA, 94550, USA  \\[3pt]
  e-mail: mgulian@sandia.edu \\[7pt]
$^3$ Department of Mathematics \\
University of Nebraska-Lincoln \\
Lincoln, NE, 68588, USA  \\[3pt]
  e-mail: hayley.olson@huskers.unl.edu \\[7pt]
$^4$ Division of Applied Mathematics, 
Brown University \\
Providence, RI, 02912, USA  \\
and Pacific Northwest National Laboratory \\
Richland, WA, 99354, USA  \\[3pt]
  e-mail: george\_karniadakis@brown.edu

\end{document}

%% file: fcaa_style_DG.tex
\usepackage{graphicx}
\usepackage{epsfig}

\usepackage{amsthm}
\usepackage{amsmath}
\usepackage{latexsym}
\usepackage{amsfonts}
\usepackage{amssymb}

\usepackage{appendix}
\usepackage{bm}
\usepackage{dsfont}
\usepackage{mathtools}
\mathtoolsset{showonlyrefs=true}
\usepackage{bbm}
\usepackage[usenames,dvipsnames]{xcolor}

\newtheoremstyle{theorem}
  {15pt}          
  {15pt}  
  {\sl}  
  {\parindent}
  {\sc}  
  {. }   
  { }    
  {}     
\theoremstyle{theorem}
\newtheorem{lemma}{Lemma}[section]
\newtheorem{theorem}{Theorem}[section]
\newtheorem{corollary}{Corollary}[section]

\newtheoremstyle{defi}
  {15pt}          
  {15pt}  
  {\rm}  
  {\parindent}     
  {\sc}  
  {. }    
  { }    
  {}     
\theoremstyle{defi}

\newtheorem{remark}{Remark}[section]

 \def\proofname{\indent\rm P\ r\ o\ o\ f.}
 \def\proofend{\hfill$\Box$}
 \def\qed{\hfill$\Box$} 
\newcommand{\Ds}{{\rm div}^s}
\newcommand{\Gs}{{\rm grad}^s}
\newcommand{\Dsl}{{\rm div}^s_\lambda}
\newcommand{\Gsl}{{\rm grad}^s_\lambda}
\newcommand{\Dsd}{{\rm div}^s_\delta}
\newcommand{\Gsd}{{\rm grad}^s_\delta}

\newcommand{\mcR}{\mathcal{R}}

\newcommand{\mcD}{\mathcal{D}}
\newcommand{\mcG}{\mathcal{G}}
\newcommand{\mcDw}{\mathcal{D}_\omega}
\newcommand{\mcGw}{\mathcal{G}_\omega}
\newcommand{\mcLw}{\mathcal{L}_\omega}

\newcommand{\mcF}{\mathcal{F}}

\newcommand{\mcL}{\mathcal{L}}

\newcommand{\mcA}{\mathcal{A}}

\newcommand{\mbR}{\mathbb{R}}

\newcommand{\mbRn}{{\mathbb{R}^n}}

\newcommand{\omg}{{\Omega}}
\newcommand{\omgi}{{\Omega_I}}

\newcommand{\omgomgi}{{\Omega\cup\Omega_I}}

\def \alphab{{\boldsymbol\alpha}}

\def \thetab{{\boldsymbol\theta}}
\def \etab {\mbox{{\boldmath$\eta$}}}

\def \wb{\mathbf{w}}
\def \vb{\mathbf{v}}
\def \xb{\mathbf{x}}
\def \yb{\mathbf{y}}
\def \zb{\mathbf{z}}

\def \xib{{\boldsymbol\xi}}

\def \etab{{\boldsymbol\eta}}

\def \rhob{{\boldsymbol\rho}}

\def \gameq{\gamma_{\textup{eq}}}

\newcommand{\vertiii}[1]{{\left\vert\kern-0.25ex\left\vert\kern-0.25ex\left\vert #1 
    \right\vert\kern-0.25ex\right\vert\kern-0.25ex\right\vert}}
    
    \newtheorem{manuallemmainner}{Lemma}
\newenvironment{manuallemma}[1]{%
  \renewcommand\themanuallemmainner{#1}%
  \manuallemmainner
}{\endmanuallemmainner}

\newtheorem{conjecture}{Conjecture}[section]